\documentclass{article}
\usepackage{amsthm,amssymb, amsmath, amscd}
\usepackage{tikz-cd}
\usepackage{cleveref}
\usepackage{enumerate}
\usepackage{todonotes}
\usepackage{mathtools}

\usepackage{mathrsfs}

\usepackage{geometry}

\makeatletter
\def\thanks#1{\protected@xdef\@thanks{\@thanks
		\protect\footnotetext{#1}}}
\makeatother

\numberwithin{equation}{section}

\theoremstyle{plain} 
\newtheorem{theorem}{Theorem}[section]
\newtheorem{lemma}[theorem]{Lemma}

\newtheorem{proposition}[theorem]{Proposition}

\theoremstyle{definition} 
\newtheorem{definition}[theorem]{Definition}

\theoremstyle{remark}
\newtheorem{remark}[theorem]{Remark}

\newcommand{\bbA}{\mathbb A}

\newcommand{\bbC}{\mathbb C}

\newcommand{\bbF}{\mathbb F}
\newcommand{\bbG}{\mathbb G}
\newcommand{\bbH}{\mathbb H}
\newcommand{\bbL}{\mathbb L}

\newcommand{\bbQ}{\mathbb Q}
\newcommand{\bbR}{\mathbb R}
\newcommand{\bbV}{\mathbb V}

\newcommand{\bbZ}{\mathbb Z}

\newcommand{\calA}{\mathcal A}

\newcommand{\calC}{\mathcal C}

\newcommand{\calL}{\mathcal L}
\newcommand{\calM}{\mathcal M}

\newcommand{\calO}{\mathcal O}
\newcommand{\calV}{\mathcal V}
\newcommand{\calW}{\mathcal W}
\newcommand{\calX}{\mathcal X}
\newcommand{\calY}{\mathcal Y}

\newcommand{\lie}[1]{ \mathfrak{#1}}
\newcommand{\psm}[1]{\left( \begin{smallmatrix}#1 \end{smallmatrix} \right)}

\newcommand{\isomto}{\stackrel{\sim}{\longrightarrow}}

\newcommand{\Ta}{\mathrm{Ta}}

\newcommand{\bad}{\mathrm{bad}}

\DeclareMathOperator{\Spec}{Spec}
\DeclareMathOperator{\Lie}{Lie}
\DeclareMathOperator{\Fil}{Fil}
\DeclareMathOperator{\Pic}{Pic}

\DeclareMathOperator{\GL}{GL}
\DeclareMathOperator{\End}{End}
\DeclareMathOperator{\Res}{Res}
\DeclareMathOperator{\Hom}{Hom}
\DeclareMathOperator{\ord}{ord}
\DeclareMathOperator{\Aut}{Aut}
\DeclareMathOperator{\Gal}{Gal}

\DeclareMathOperator{\GSpin}{GSpin}
\DeclareMathOperator{\deghat}{\widehat{\mathrm{deg}}}

\newcommand{\inv}{\mathrm{inv}}
\newcommand{\alg}{\mathrm{alg}}

\title{Zero-dimensional Shimura varieties and central derivatives of Eisenstein series }
\author{Siddarth Sankaran\thanks{This research was supported by an NSERC Discovery grant.}}
\date{}
\begin{document}
		\maketitle
	
	\begin{abstract}
		We formulate and prove a version of the arithmetic Siegel-Weil formula for (zero dimensional) Shimura varieties attached to tori, equipped with some additional data. More precisely, we define a family of ``special"  divisors in terms of Green functions at archimedean and non-archimedean places,  and prove that their degrees coincide with the Fourier coefficients of the central derivative of an Eisenstein series. The proof relies on the usual Siegel-Weil formula to provide  a direct link between both sides of the identity, and in some sense, offers a  more conceptual point of view on prior results in the literature.
	\end{abstract}

	\tableofcontents
	
{\ \\}

	\section{Introduction}
			In this paper, we formulate and prove a version of the ``arithmetic Siegel-Weil formula" for zero-dimensional Shimura varieties. Broadly speaking, the arithmetic Siegel-Weil formula, which is a conjectural formula due to Kudla, predicts that the arithmetic heights of certain ``special" cycles on  Shimura varieties are related to the Fourier coefficients of derivatives of Eisenstein series. For zero-dimensional varieties, a number of special cases have been worked out in the literature; see for example \cite{KRY-IMRN, HowardCM,HowardYang} as well as \cite{AGHMP}, where an identity of this form was a crucial ingredient in their proof of the average Colmez conjecture.

		Here we take a more general, abstract point of view. Our starting point is an arbitrary Shimura datum $(T, h)$, where $T$ is a rational torus, and $h \colon \Res_{\bbC/\bbR} \bbG_m \to T_{/{\bbR}}$ is a morphism. Upon fixing a compact open subgroup $K \subset T(\bbA_f)$, this data determines a zero-dimensional Shimura variety that admits a canonical model over the reflex field $E(T, h)$. Choose a CM field $E$ containing $E(T,h)$ and let $\calM$ denote the base change of the canonical model to $\Spec (E)$; then the theory of complex multiplication provides for an action of $E^{\times} \backslash \bbA_E^{\times}$ on the geometric points $\calM(E_v^{\alg})$ for any place $v$ of $E$. 
		
		We also fix the following collection of data, which we describe in loose terms here, and refer to \Cref{sec:data} for more precise details: 
			\begin{itemize}
				\item   an ``{incoherent family}" $\calV = (\calV_v)_v$ whose members are one-dimensional local $E_v$-Hermitian spaces, one for each place $v$ of $F$ (here $F \subset E$ is the maximal totally real subfield); 
				\item a family of lattices $\calL_v \subset \calV_v$ for non-archimedean places $v$; and
				\item a family of pairs $\left\{(\bbV^{(v)}, \beta^{(v)} )\right\}$ indexed by the places of $F$ that are non-split in $E$; here $\bbV^{(v)}$ is a local system of one-dimensional $E$-Hermitian spaces on the base change $\calM_{/ E_{v}^{\alg}}$, and   $\beta^{(v)}$ is a map relating this local system to the incoherent collection $\calV$.
				In particular, at each geometric point $z \in \calM(E_v^{\alg})$, the map $\beta^{(v)}$ identifies the fibre of $\bbV^{(v)} \otimes \bbA_E$ at $z$ with the coherent collection that differs from  $\calV$ at precisely the place $v$. 
			\end{itemize}

		In addition,
		the action of $E^{\times} \backslash \bbA_E^{\times}$ on the geometric points of $\calM$ restricts to an action of $H(F) \backslash H(\bbA_F)$, where, for an $F$-algebra $R$, we write $H(R) = \{ x \in  R \otimes_F E \mid N_{E/F}(x) = 1  \}$; we further require that the data above satisfies a certain equivariance property with respect to the action of $H(F) \backslash H(\bbA_F)$.

		With this data in place, and for a fixed $\alpha \in F^{\times}$, we construct a family of ``Green functions"
			\begin{equation}
				g_v(\alpha, \tau) \colon \calM(E_v^{\alg}) \to \bbR
			\end{equation} 
			indexed by the places $v$ of $F$, and depending on an additional parameter $\tau \in \bbH_E$, the Hilbert upper-half space attached to $E$. When $v$ is an archimedean place, this coincides with the Green function defined by Kudla. The novelty of the present article is an analogous construction at  non-archimedean places, which can be roughly described as a non-archimedean factor of an incomplete Mellin transform.
			Packaging these Green functions together as $v$ varies over all the places of $F$ yields a ``special cycle" 
				\begin{equation}
					\widehat Z(\alpha, \tau) \in \widehat \Pic_{\bbA}(\calM)
				\end{equation}
		 where $\widehat \Pic_{\bbA}(\calM)$ is the group of isomorphism classes of adelically metrized line bundles.

			On the other hand, we use the same data to define a Hilbert modular Eisenstein series $E(\tau, s, \Phi^*)$ of parallel weight one that vanishes at the centre $s=0$ of its functional equation. Our main theorem relates the two constructions:
			
				\begin{theorem} 
						For every $\alpha \in F^{\times}$, we have
							\begin{equation}
								\widehat {\mathrm{deg}} \,  \widehat Z(\alpha, \tau) = -   \# \calM(\bbC)  \cdot E'_{\alpha}(\tau, 0, \Phi^*)
							\end{equation}
						where $E'_{\alpha}(\tau, 0, \Phi)$ is the $\alpha$'th Fourier coefficient of the central derivative
							\begin{equation}
								E'(\tau, 0, \Phi^*) := \frac{d}{ds} E(\tau, s, \Phi^*)\big|_{s=0}.
							\end{equation}
				\end{theorem}

		The  aforementioned special cases can be placed in the context of our main theorem (though the translation is not entirely trivial), and in particular we obtain new proofs of those results. It should be emphasized that our approach circumvents the explicit deformation-theoretic calculations that figure prominently in the prior literature. Instead, in these cases, the $p$-adic Green functions  we construct encode arithmetic degrees in a natural way (see \Cref{prop:Green Arakelov}), and our main identity is proved via the Siegel-Weil formula, without requiring an explicit computation of the left hand side. In particular, this line of argument provides a more conceptual explanation for such identities. We present a few examples illustrating this point in the last section of the paper.

	\section{Preliminaries}
	 
	\subsection{Adelic line bundles on zero-dimensional spaces and arithmetic degrees} \label{sec:adelic metric}
		In this section, we give a very naive description of the theory of adelic line bundles, as introduced by Zhang, see e.g.\ \cite{ZhangSmallPoints}: as our interest lies solely in the zero-dimensional case, we will give an ad-hoc presentation here relevant to this setting. 
		
		Suppose $X$ is a zero dimensional variety over a number field $L$, and $\omega$ is a line bundle on $X$. 
		For a place $v$ of $L$, choose an algebraic close $L_v^{\alg}$ of $L_v$. A \emph{v-adic metric} 
		\begin{equation} 
		\| \cdot \|_v = \{ \| \cdot \|_x \}_{x \in X( L_v^{\alg})}
		\end{equation} 
		on $L$ is  a collection of (real-valued) metrics
		\begin{equation}
		\| \cdot \|_{x} \colon \omega_x \to \bbR_{\geq 0}, \qquad x \in X( L_v^{\alg})
		\end{equation}
		on the fibres $\omega_x$ of $\omega$ at geometric points. These metrics are required to satisfy 
		\begin{equation}
		\| \lambda s \|_x =  |\lambda|_v  \cdot \| s \|_x ,  
		\end{equation}
		for all $ \lambda \in L_v^{\alg}, s \in \omega_x$, where $|\cdot |_v$ is the unique extension of the norm on $L_v$ to $L_v^{\alg}$.

		At a non-archimedean place $v$, a key role is played by \emph{model metrics}, defined as follows. Suppose $X_v =  X \times_L \Spec(L_v)$ admits a projective model $\calX $ over $\Spec(\calO_{L,v})$ and that $\omega$ extends to a line bundle $\tilde \omega$ over $\calX$. Then  an $L_v^{\alg}$ valued point $z  \in X(L_v^{\alg})$ lifts uniquely to a point $\tilde z \in \calX(\calO_{L_v^{\alg}})$, and hence we obtain an $\calO_{L_v^{\alg}}$-lattice
		\begin{equation}
		\widetilde{\omega}_{\tilde z} \subset  \omega_z.
		\end{equation}
		This inclusion determines a metric, called the \emph{model metric relative to $(\calX, \tilde \omega)$}, by the formula
		\begin{equation}
		\| s \|_{z, (\calX,\tilde \omega)} := \inf \{ |a|_v \ | \ a \in (L_{v}^{\alg})^{\times}, \ a^{-1} s \in \widetilde\omega_{\tilde z}    \}
		\end{equation}

		\begin{definition} An  adelic line bundle $\widehat \omega = ( \omega, (\| \cdot \|_v)_v)$ consists of a line bundle $\omega$ together with $v$-adic metrics $\| \cdot \|_v$ for every place $v$, such that the following condition holds: there exists an integer $N$ and a model $(\calX, \tilde \omega)$ over $\Spec \,  \calO_L[1/N]$ such that the metric $\| \cdot \|_v $ is induced by this model at every non-archimidean place induced by a prime ideal relatively prime to $N$. 
			
			Two adelic line bundles $\widehat \omega $ and $\widehat \omega' $ are \emph{isomorphic} if there is an isomorphism $\omega \simeq \omega'$ that induces an isometry $ (\omega_x, \| \cdot \|_x) \simeq (\omega'_x, \|\cdot\|'_x)$ at every geometric point $x$.  
			
		\end{definition}
	
			Of course, for zero dimensional schemes these definitions are overkill, as all line bundles are isomorphic to the trivial one. When $\omega = \calO_X$ is the trivial bundle, specifying an adelic metric $(\| \cdot \|_v)_v$ is equivalent to choosing a family $(g_v)$ of ``Green functions" $g_v \colon X(L_v^{\alg}) \to \bbR$ whose members are identically zero for almost all $v$; this equivalence is encoded by the relation
			\begin{equation}
				\| 1 \|^2_z = e^{- g_v(z)}, \qquad z \in X(L_v^{alg}).
			\end{equation}

			Let $\widehat\Pic_{\bbA}(X)$ denote the group of isomorphism classes of  adelic line bundles; as every bundle on $X$ is isomorphic to the trivial one, every element of $\widehat\Pic_{\bbA}(X)$ is determined uniquely by a set of Green functions $g_v \colon X(L_v^{\alg}) \to \bbR$ as above.
			
			Given an adelic bundle $\widehat \omega $, we define its arithmetic degree
				\begin{equation}
					\widehat \deg \, \widehat \omega = \sum_{v \leq \infty} \sum_{z \in X(L_v^{\alg})} - \log \| \eta \|_z^2
				\end{equation}
			where $\eta$ is a non-vanishing section of $\omega$. The product formula implies that this map is invariant on isomorphism classes, and therefore defines a function on $\widehat\Pic_{\bbA}(X)$. In particular, if $\widehat \omega = (\calO_X, (g_v)_v) \in \widehat \Pic_{\bbA}(X)$ is given by a collection of Green functions   as above, then
				\begin{equation}
					\widehat \deg \, \widehat \omega = \sum_{v \leq \infty} \sum_{z \in X(L_v^{\alg})} g_v(z).
				\end{equation}
				
	\subsection{Hermitian spaces} \label{sec:prelim Hermitian}
 	Let $E$ be a CM field, with maximal totally real subfield $F$. Suppose that $v$ is a place of $F$, let $E_v := E \otimes_F F_v$ and suppose that $V_v$ is a Hermitian space of dimension $m$ over $E_v$. We define the local invariant
	\begin{equation}
	\inv_{v}(V_v) := ( (-1)^{m(m-1)/2} \det V_v, \Delta_{E/F})_{v} \in F_v^{\times}/ N_{E_v/F_v}(E_v^{\times})
	\end{equation}
	where $\Delta_{E/F}$ is the discriminant, and $(\cdot, \cdot)_v$ is the Hilbert symbol. At a non-archimedean place, two $E_v$-Hermitian spaces of the same dimension  are isomorphic if and only if they have the same invariants. In particular, at split places, there is only one isometry class of $E_v$-Hermitian space of a given dimension, while at non-split places, there are two. At archimedean places, the isometry class of a Hermitian space over $E_v \simeq \bbC$ is determined by its signature. 
	
	Now suppose $V$ is a global Hermitian space over $E$. Almost all of its local invariants are equal to one, and together they satisfy the product formula 
	\begin{equation}	
	\prod_{v \leq \infty}  \inv_v(V) = 1.
	\end{equation}
	Conversely, suppose we have a collection of signs  $\epsilon_v \in \{ \pm 1 \}$ for $v \leq \infty$, almost all of which are 1, and a collection of pairs of non-negative integers $(r_v, s_v)$ for $v | \infty$, with $r_v + s_v = n$ independent of $v$. Then there exists  a Hermitian space $V$ of signature $(r_v, s_v)$ at each $v | \infty$ and $\inv_v(V) = \epsilon_v$ if and only if $\prod_{v \leq \infty} \epsilon_v = 1$ and $\epsilon_v = (-1)^{s_v+1}$ for each $v |\infty$. 
	
	Following the terminology of Kudla  \cite{KudlaCentralDerivs}, we say that a collection
		\begin{equation}
			\calV = (\calV_v)_{v \leq \infty}
		\end{equation}
	of local Hermitian spaces is \emph{incoherent} if there is no global Hermitian space that is locally isometric to $\calV_v$  for all $v$.
	
	\subsection{Whittaker functions and Eisenstein series}
We briefly recall the theory of (Hermitian) Hilbert Eisenstein series, mainly to fix notation.   We fix, once and for all, a character
	\begin{equation}
		\chi \colon  E^{\times} \backslash \bbA_E^{\times} \to \bbC^{\times}.
	\end{equation}

Consider the group $G = U(1,1)$, viewed as an algebraic group over $F$. Concretely, for an $F$-algebra $R$, the points of $G$ are
	\begin{equation}
		G(R) \ := \ \left\{  g \in \GL_2(R \otimes_F E) \ | \ g \psm{ & 1 \\ 1 & } = \psm{ & 1 \\ 1 & } \overline g^t  \right\}
	\end{equation}
Let $B \subset G$ denote the subgroup of upper-triangular matrices, with Levi decomposition $B = MN$, where, for an $F$-algebra $R$, we have
	\begin{equation}
		M(R) = \left\{ m(a) = \psm{a & \\ & \overline a^{-1}} \ | \ a \in (R \otimes_F E)^{\times} \right\}
	\end{equation}
and
	\begin{equation}
		N(R) = \left\{ n(b) = \psm{1 & b \\ & 1 } \ | \ b \in R  \right\}.
	\end{equation}

Now fix a place $v$ of $F$. For every $s \in \bbC$, we obtain a character
\begin{equation}
\chi | \cdot |^s \colon B(F_v) \to \bbC^{\times}, \qquad \psm{a & b \\ & \overline a^{-1}} \mapsto \chi_v(a) | a|_{v}^s, 
\end{equation}
	where  $a \in E_v, \, b \in F_v$.
We define the degenerate principal series representation to be the (smooth) normalized induction
\begin{equation}
I_v(s,\chi)  := \mathrm{Ind}_{B_v}^{G_v} \left(  \chi_v | \cdot |_v^s \right)
\end{equation}
where we abbreviate $G_v = G(F_v)$ and $B_v = B(F_v)$.
Concretely, elements of $I_v(s,\chi)$ are smooth function $\Phi(\cdot  , s) \colon G_v \to \bbC$ such that
\begin{equation} \label{eqn:deg ps trans law}
\Phi\left( \psm{a & b \\ & a^{-1}} g, s \right) \ = \ \chi_v(a) | a|_v^{s+1} \, \Phi(g, s),
\end{equation}
for all  $\psm{a & b \\ & a^{-1}} \in B_v$.

We say that a section $\Phi(g,s)$ is holomorphic (resp.\ meromorphic) if for every fixed $g$, the map $s \mapsto \Phi(g,s)$ is holomorphic (resp.\ meromorphic). Additionally, consider the maximal compact subgroup\footnote{In the archimedean case, we view $U(1) \times U(1)$ as a subgroup of $G(\bbR)$ via the map 
		\[
			(e^{i\theta}, e^{i \varphi}) \mapsto e^{i \varphi} \psm{ \cos \theta & \sin \theta \\ - \sin \theta & \cos \theta }
		\].
}
	\begin{equation}
		K_v = 
			\begin{cases} 
				G(\calO_{F,v}), & v \text{ non-archimedean,} \\ 
				U(1) \times U(1) , &  v \text{ archimedean};
			\end{cases}	
	\end{equation}
we say $\Phi(g,s)$ is \emph{standard} if $\Phi(k,s)$ is independent of $s$ for all $k \in K_v$. 

An important source of standard sections is the Weil representation: suppose that $V_v$ is a Hermitian space over $E_v$ of dimension $m$, and assume that the restriction of the character $\chi$ to $\bbA_F$ coincides with $(\chi_{E/F})^m$, where $\chi_{E/F} \colon F^{\times} \backslash \bbA_F^{\times} \to \{ \pm 1\}$ is the quadratic character corresponding to the extension $E/F$.

Let  $S(V_v)$ denote the space of Schwartz-Bruhat functions on $V_v$, and let 
	\begin{equation}
		\omega_{\calV_v} \colon G(F_v) \to \Aut S(V_v)
	\end{equation} denote the Weil representation, cf.\  \cite{KudlaSplitting}. This representation depends on the additional choice of an additive character $\psi$ on $\bbA_F$; we take $\psi = \psi_{\bbQ} \circ \mathrm{tr}_{F/\bbQ}$, where $\psi_{\bbQ}$ is the standard additive character on $\bbA_{\bbQ}$, and suppress this choice from the notation.

Explicit formulas for the action of elements of $M_v$ and $N_v$ are given as follows:
	\begin{align}
		\omega_{\calV_v}\left( m(a) \right)  \varphi_v(x) \ &= \ |a|_v^{m} \chi_v( a)  \varphi_v(ax), & a &\in E_v^{\times} \\
		\omega_{\calV_v}\left( n(b) \right)  \varphi_v(x) \ &  = \  \psi_{v}(bQ(x)) \varphi_v(x), & b & \in F_v.
	\end{align}
Comparing these formulas with \eqref{eqn:deg ps trans law}, we find that there is an $G_v$-equivariant map
	\begin{equation}
		\lambda_{V_v} \colon S(V_v) \to I_v(s_0, \chi), \qquad s_0 := m-1
	\end{equation}
 given by the formula
	\begin{equation}
		\lambda_{V_{v}} (\varphi )(g) =  \left(\omega_{V_v}(g) \varphi \right)(0).
	\end{equation}
Given $\varphi \in S(V_v)$, we denote by
	\begin{equation}
		\Phi(g, s, \lambda_{V_v}(\varphi)) \in I_v(s, \chi)
	\end{equation}
the unique standard section whose value at $s_0$ coincides with $\lambda_{V_v}(\varphi)$. 

Next, for $m \in F_v$ and $\Phi(g,s) \in I_v(s,\chi)$, consider the Whittaker functional
	\begin{equation}
		W_{m,v}(g, s, \Phi) := \int_{F_v} \Phi(w n(b)g, s) \, \psi_v(-mb) \, db
	\end{equation}
where $w = \psm{ & -1 \\  1 & }$ and the measure $db$ is normalized so that it is self-dual for the pairing $(b_1, b_2) \mapsto \psi_v(b_1b_2)$.

\begin{proposition} \cite{KudlaSweet} \label{prop:local whittaker dichotomoy}
	Suppose $v$ is a non-split place of $E$, and let $V_v$ be a one-dimensional Hermitian space over $E_v$. Let $R(V_v)$ denote the image of the Rallis map $\lambda_{V_v}$. 
	\begin{enumerate}[(i)]
		\item Let $m \in F_v$. If $V_v$ does not represent $m$, then 
			\begin{equation}
				W_{m,v}(g, 0, \Phi) =0
			\end{equation}
		for all $\Phi \in R(V_v)$. 
		\item Let $V^+_v$ and $V^-_v$ denote two one-dimensional Hermitian spaces with $\inv_v (V^{\pm}_v) = \pm 1$. Then
			\begin{equation}
				I_v(0, \chi) = R(V_v^+) \oplus R(V_v^-).
			\end{equation} 
	\end{enumerate}
\qed
\end{proposition}

We now turn to the global setting. Let
	\begin{equation}
			I(s, \chi) = \mathrm{Ind}_{B(\bbA_F)}^{G(\bbA_F)} \left( \chi \, | \cdot |_{\bbA}^s \right)  
	\end{equation}
denote the global degenerate principal series representation consisting of smooth functions $\Phi(g,s) \colon G(\bbA_F) \to \bbC$ satisfying
	\begin{equation}
		 \Phi \left( \psm{a & b \\ & \overline a^{-1}} g, s \right) = \chi(a) |a|^{s+1} \Phi(g,s) 
	\end{equation}
for all $ \psm{a & b \\ & \overline a^{-1}} \in B(\bbA_F)$. As before, a section $\Phi(g,s) \in I(s,\chi)$ is said to be holomorphic if for each fixed $g$, it defines a holomorphic function in $S$, and is standard if the restriction to $K = \prod_v K_v$ is independent of $s$.  There is a natural factorization
	\begin{equation}
		I(s, \chi) = \bigotimes_v I_v(s, \chi)
	\end{equation}
and we say a section $\Phi(g,s)$ is \emph{factorizable} if it can be written as a product $\Phi(g,v) = \otimes_v \Phi_v(g,s)$.

If $\Phi(g,s) \in I(s,\chi)$ is a holomorphic standard section, we define the Eisenstein series
	\begin{equation}
		E(g, s, \Phi) := \sum_{\gamma \in B(F) \backslash G(F) }  \Phi(\gamma g, s).
	\end{equation}
Standard results in the theory assert that this series is absolutely convergent for $Re(s) \gg 0$, has meromorphic continuation to $s \in \bbC$, and  a functional equation $s \leftrightarrow -s$; in particular, our Eisenstein series have been normalized so their centre of symmetry is at $s=0$. 

Writing the Fourier expansion of $E(g, s, \Phi)$ as
	\begin{equation}
		E(g,s,\Phi) = \sum_{m \in F} E_m(g,s,\Phi),
	\end{equation}
and assuming that $\Phi(g,s)$ is factorizable, we have that
	\begin{equation} \label{eqn:Eis product}
		E_m(g, s, \Phi) = W_{m}(g, s, \Phi) = \prod_v W_{m,v}(g_v, s, \Phi_v)
	\end{equation}
for all $m \neq 0$. 

We conclude this section with a discussion of the values of Eisenstein series at the central point $s=0$, in the cases of relevance to this paper.

Suppose that $\calV = (\calV_v)_v$ is a collection of one-dimensional local Hermitian spaces indexed by the places of $F$, and we are given a Schwartz function
	\begin{equation}
		\varphi = \otimes \varphi_v \in S(\calV), \qquad \varphi_v \in S(\calV_v)
	\end{equation}
with corresponding standard section $\Phi(g, s, \lambda_{\calV}(\varphi)) = \otimes_v \Phi_v(g, s, \lambda_{\calV_v}(\varphi_v))$.

If $\calV$ is an \emph{incoherent} collection, then for any $m \in F^{\times}$, there is at least one place $v$ such that $\calV_v$ does not represent $m$; in light of \eqref{eqn:Eis product} and \Cref{prop:local whittaker dichotomoy}(i), we find
	\begin{equation}
		E(g, 0, \Phi) = 0
	\end{equation}
	
When $\calV$ is a coherent collection, the central value is computed by the {Siegel-Weil} formula, as follows. Suppose $V$ is a global one-dimensional Hermitian space over $E$. Let $H$ denote the algebraic group over $ F$ whose points, for an $F$-algebra $A$, are given by
	\begin{equation} \label{eqn:def of H}
		H(A) = \{x \in A \otimes_F E | \  \ x \overline x  = 1 \};
	\end{equation} 
in particular, we may identify $H(A)$ as the unitary group of $V \otimes_F A$. 

For an Schwartz function $\varphi \in S(V\otimes_F \bbA_F)$ and an element $h \in H(\bbA_F)$, we define the theta series
	\begin{equation} \label{eqn:theta function defn}
		\Theta(g, h, \varphi) = \sum_{x \in V} \omega_V(g) \varphi(h^{-1} x)
	\end{equation}
which defines an automorphic form on $G \times H$. 

The following proposition is a special case of the Siegel-Weil formula, which in the unitary setting is due to Ichino.
\begin{proposition}[{\cite[Theorem 4.2]{Ichino2004}}] \label{prop:SW}  
	Let $V$ be a one-dimensional Hermitian space over $E$, let $\varphi \in S(V\otimes \bbA_E)$, and let $\Phi(g,s) = \Phi(g,s, \lambda_V(\varphi))$ denote the corresponding section. Then
		\begin{equation}
			E(g, 0, \Phi) = 2 \int_{[H]} \, \Theta(g, h, \varphi) \, dh
		\end{equation}
	where $dh$ is the Haar measure on $[H] = H(F) \backslash H(\bbA_F)$ normalized so that the total volume is one. \qed
\end{proposition}

	\section{Description of the data} \label{sec:data}
	In this section, we give a detailed description of the data and assumptions that go into our main theorem. 
	
	\subsection{Zero dimensional Shimura varieties}
			Let $T$ be a rational torus, and $h \colon \mathbb S \to (T)_{\bbR} $ a homomorphism, where $\mathbb S = \Res_{\bbC/\bbR} \bbG_m$ is Deligne's torus.
			
			For a neat compact open subgroup $K \subset T(\bbA_{\bbQ,f})$, we have a zero dimensional Shimura variety
				\begin{equation}
					Sh_K(T, h) := T(\bbQ) \big\backslash \{ h\} \times T(\bbA_{\bbQ, f}) / K.
				\end{equation}
Suppose that $E \subset \bbC$ is a CM field that contains the reflex field  $E(T,h)$, so that $Sh_K(T,h)$ admits a canonical model $\calM$ over $\Spec(E)$. As per the theory of canonical models, there is an action of $E^{\times} \backslash \bbA_E^{\times}$  on $\calM(\bbC) = \calM(E^{\alg})= Sh_K(T,h)$ that factors through the Galois group $\Gal(E^{\alg}/E)^{\mathrm{ab}}$, where $E^{\alg}$ is a fixed algebraic closure of $E$, see e.g.\ \cite{MilneNotes} for details. For any place $v$ of $E$, we fix an embedding $E^{\alg} \subset E^{\alg}_v$, so  we may view $E^{\times} \backslash \bbA_E^{\times}$ acting on $\calM(E_v^{\alg}) = \calM(E^{\alg})$  as well. 
			
			Let $F \subset E$ denote the maximal totally real subfield, and consider the algebraic group $H$ over $F$  defined by \eqref{eqn:def of H}; viewing $H(\bbA_F)$ and $H(F)$ as a subgroups of $\bbA_E^{\times}$ and $E^{\times}$, respectively, we have a natural inclusion $H(F) \backslash H(\bbA_F) \hookrightarrow E^{\times} \backslash \bbA_E^{\times}$, and so in this way, we obtain an action of $H(F) \backslash H(\bbA_F) $ on $\calM(E_v^{\alg})$.
			
			
	\subsection{Linear algebraic data} \label{sec:lin alg}
			Let $v$ be a place of $F$, and set $E_v := E \otimes_F F_v$. We say that the place $v$ is \emph{non-split} if $E_v$ is a field; in particular, the archimedean places of $F$ are all non-split.
			
			Our next piece of data is an \emph{incoherent collection} 
				\begin{equation}
					\calV  = (\calV_v)_{v \leq \infty}
				\end{equation}
			indexed by the places of $F$; each member  $(\calV_v, \langle \cdot, \cdot \rangle_v)$ is a Hermitian space over $E_v$, and we assume  that 	
				\begin{equation}
					\calV_v \text{ is positive definite  for all } v | \infty.
				\end{equation}

			For each non-archimedean place $v < \infty$, we also fix an $\calO_E$-stable lattice
				\begin{equation}
					\calL_v \subset \calV_v
				\end{equation}
			and assume that $\calL_v$ is  self-dual  for almost all $v$. Note that we are not assuming that $\calL_v$ is integral for all $v$.  
			
			For each non-archimedean place $v$, we define the Schwartz function 
				\begin{equation}
					\varphi_{\calL,v} \in S(\calV_v)
				\end{equation}
			to be the characteristic function of $\calL_v$. At an archimedean place $v$, we set 
				\begin{equation}
					\varphi_{\calV,v}(x) = e^{- \pi \langle x, x \rangle_v}.
				\end{equation}

			For convenience, we set
				\begin{align} \label{eqn:Schwartz L defn}
					\varphi_{\calL,f} &= \bigotimes_{v < \infty} \varphi_{\calL,v},	&  
						\varphi_{\calV,{\infty}} &= \bigotimes_{v | \infty} \varphi_{\calV,v},	 &
					 \varphi_{\calL,\bbA} &=  \varphi_{\calV,\infty} \otimes \varphi_{\calL,f}.
				\end{align}
			 
			  Next, let $v$ be a  non-split place of $F$. Define an element $\varsigma_{v} \in F_{v}^{\times}$ as follows: if $v$ is non-archimedean and inert, let $\varsigma_{v} = \varpi$, where $\varpi$ is a uniformizer of $F_{v}$. If $v$ is ramified, choose $\varsigma_{v} \in \calO_{F, v}^{\times}$ such that $\varsigma_{v} \notin N(E^{ \times}_{v})$. If $v$ is archimedean, we  take $\varsigma_v = -1$.
			 
			 In any of the above cases, let $\calV'_{v}$ denote the $E_v$-Hermitian space whose underlying space is $\calV_v$, with Hermitian form  given by 
			 \begin{equation}
					\langle x, y \rangle'_v := \varsigma_v \cdot \langle x, y \rangle_v
			 \end{equation}
			 By construction, $\calV_{v}$ and $\calV'_{v}$ are representatives for the two distinct isomorphism classes of $E_{v}$ Hermitian spaces of dimension 1. 
			
			Finally, we note that as $\calV_v$ and $\calV'_v$ share the same underlying vector space, we view $\calL_v$ as a lattice in either space; in particular, we may  view $\varphi_{\calL,v}$ as a Schwartz function for $\calV'_v$ as well.
			
		\subsection{Local systems}
			 For each non-split  place $v$ of $F$, fix an algebraically closed extension $E_{v}^{\alg}$ of $E_v$, and  set 
				\begin{equation}
					\calM_{/E_{v}^{\alg}} := \calM \times_E \Spec(E_{v}^{\alg}).
				\end{equation}

			Our final piece of data is a family of pairs
				\begin{equation}
				 \left\{ (  \bbV^{(v)}, \,  \beta^{(v)}) \ | \ v \text{ non-split } \right \}
				\end{equation}
			 	whose entries are as follows:
					\begin{enumerate}[(a)]
						\item $\bbV^{(v)}$ is a   local system of one-dimensional $E$-Hermitian spaces  over $\calM_{/E_{v}^{\alg}}$  equipped with an $H(\bbA_F)$ action lifting the action on $\calM_{/E_{v}^{\alg}}$.
						\item Let 
							\begin{equation}
						   		\calV^{v} \times \calV'_{v} = \left(\prod_{\substack{ w \leq \infty \\ w \neq v}} \calV_w \right)  \times \calV'_{v}
							\end{equation}  
						which we view as a Hermitian $\bbA_E $-module, and denote by $\underline{\calV^{v} \times \calV'_{v} }$  the constant local  system on $	\calM_{/E_{v}^{\alg}}$, whose fibre at each point is a copy of 	$\calV^{v} \times V'_{v}$
						
						 The datum $ \beta^{(v)}$ is an  isometry
						 	\begin{equation}
							  \beta^{(v)} \colon \,   \bbV^{(v)} \otimes_E \bbA_E \to \underline{ \calV^{v} \times \calV'_{v}}
						 	\end{equation}
						 	of Hermitian $\bbA_E$-local systems, for which we require the following equivariance property to hold: given $z \in \calM(E_{v}^{\alg})$ and $t \in H(\bbA_F)$, we have an isometry 
						 	\begin{equation} 
						 		r(t) \colon 	 \bbV^{(v)}|_z   \to 	\bbV^{(v)}|_{ t \cdot z}   
						 	\end{equation}
						 	arising from part (a) of the definition, and we require that the diagram
						 		\begin{equation} \label{eqn:beta equivariance}
							 		\begin{CD}
							 			\bbV^{(v)}|_z \otimes_E \bbA_E @>  \beta^{(v)}_z >> \calV^{v} \times \calV'_{v} \\
							 			@Vr(t)\otimes 1VV @VV \times t V \\
							 			 \bbV^{(v)}|_{ t \cdot z} \otimes_E \bbA_E @> 	\beta^{(v)}_{ t \cdot z} >> \calV^{v} \times \calV'_{v}
							 		\end{CD}
						 		\end{equation}
						 	commutes up to multiplication by an element of $\widehat \calO{}^{\times}_E$. In other words, for any $t \in H(\bbA_F)$, we require that the composition
						 		\begin{equation}
							 		\beta_z^{(v)} \circ \left( r(t^{-1}) \otimes t  \right) \circ \left(  \beta^{(v)}_{t \cdot z} \right)^{-1} \ \in \ \End\left( {\calV^{v} \times \calV'_{v} } \right) 
						 		\end{equation} 
						 	is given by multiplication by an element of $\widehat \calO{}^{\times}_E$. 

					\end{enumerate}

\section{Eisenstein series and Green functions}
	We fix a collection of data 
		\begin{equation}
			\left( \calM \to \Spec (E),\,  \, \calV,   \,  (\calL_v) , \,  \left(  \bbV^{(v)},   \beta^{(v)} \right)_{v}  \right)
		\end{equation}
		as described in \Cref{sec:data}. In the following three subsections, we show how this data defines a section $\Phi^*(g,s) \in I(s,\chi)$, a collection of ``special cycles" $\widehat Z(\alpha,\tau) \in \widehat\Pic_{\bbA}(\calM_K)$, and then state and prove our main theorem.

	\subsection{The section $\Phi^*(s)$} \label{sec:Phi*}
		We begin with the non-archimedean setting. Fix a non-archimedean place $v$ of $F$ that is non-split in $E$, and let $E_v = E \otimes_F F_v$, which is a field. Then $\calV_{v}$ and $\calV'_{v}$ are the two non-isomorphic one dimensional   Hermitian spaces over $E_{v}$. 
		
		Suppose for the moment that $\calV_{v}$ represents $1$, i.e.\ there exists a vector $x \in \calV_{v}$ such that $\langle x, x \rangle_{\calV_{v}} = 1$. Then there are isometries
			\begin{equation}
				\calV_{v} \isomto \calV_{v}^+ := (E_{v}, x \overline y) \qquad \text{ and } \qquad \calV'_{v} \isomto \calV^-_{v} := (E_{v}, \varsigma_{v}x \overline y)
			\end{equation}
			where for $a \in F_{v}^{\times}$, the space $(E_{v}, a x \overline y)$  denotes the Hermitian space with underlying space $E_{v}$ and Hermitian form $\langle x, y \rangle =a x \overline y$. 
	
		If $\calV_{v}$ does not represent 1, then $\calV'_{v}$ does, and so there are isometries
			\begin{equation}
				\calV_{v} \isomto \calV^-_{v}, \qquad \calV_{v}' \isomto \calV^+_{v}.
			\end{equation}
			
		In either case, these isomorphisms carry  the lattice $\calL_{v} \subset \calV_{v}= \calV'_{v}$ to some lattice $\Lambda \subset E_{v}$. It follows easily from definitions that 
			\begin{equation}
				\Phi_{v}(g, s, \lambda_{\calV_{v}}(\varphi_{\calL_v}))  = \Phi_{v}(g, s,  \lambda_{\calV^{\epsilon}_{v}}(\varphi_{\Lambda})).
			\end{equation}
		where
			\begin{equation}
				\epsilon = \chi_{v}(\det \calV_{v}) = \pm 1. 
			\end{equation}
	Similarly, we have
		\begin{equation} \label{eqn:whittaker prime to -eps}
			\Phi_{v}\left(g, s, \lambda_{\calV'_{v}}(\varphi_{\calL_v})\right) = \Phi_{v}\left(g, s, \lambda_{\calV^{-\epsilon}_{v}}(\varphi_{\Lambda})\right)
		\end{equation}
		
		Now let $\varphi_{v}^0(x) $ denote the characteristic function of $\calO_{E,v}$. Write $\Lambda= t \cdot \calO_{E,v}$ for some $t \in E_{v}^{\times}$, so that 
			\begin{equation} \label{eqn:rescale Schwartz fn}
				\varphi_{\Lambda}(x) = \varphi_{v}^0( t^{-1} x) = |t|_{v} \, \chi(t)  \left(  \omega_{\calV^{\pm}_{v}}\left(m(t^{-1}) \right) \varphi_{v}^0 \right) (x).
			\end{equation}
		This in turn implies that   the difference
			\begin{equation}
			\delta_{v}(g, s) := \Phi_{v}(g, s, \lambda_{\calV_{v}}(\varphi_{\calL})) - |t|^{1-s}_{v} \chi(t) \Phi_v\left(g m(t)^{-1}, s, \lambda_{\calV^{\epsilon}_{v}} (\varphi^0_{v}) \right)
			\end{equation}
		vanishes identically at $s=0$. 
			\begin{remark} \label{rmk:modified local section unram}
				If $v$ corresponds to a prime of $F$ that is unramified in $E$ and ${\calL}_{v} = \calO_{E, v}$ is self-dual, then we may take $t = 1$ above, and $\delta_{v}(g,s ) \equiv 0$. In particular, $\delta_{v}(g, s) = 0$ for almost all $v$. In this case,  $\Phi_{v}(g, s, \lambda_{V_{v}}(\varphi_{\calL}))  $ equals  $\Phi_{v}(g, s, \lambda_{V^+_{v}}(\varphi^0_{v}))$, which is the standard spherical section.
			\end{remark}
			
		Next, recall that
			\begin{equation}
				I_{v}(0, \chi) = R(\calV_{v}) \oplus R(\calV_{v}'),
			\end{equation}
		so there are Schwartz functions $\phi_{v} , \phi_{v}' \in S(E_{v}) = S(\calV_{v}) = S(\calV'_{v})$   such that 
			\begin{equation}
				\delta_{v}(g,s)= s \cdot  \left(   \Phi \left(g, s, \lambda_{\calV_{v}}(  \phi _{v})\right) +   \Phi\left(g, s, \lambda_{\calV'_{v}}(  \phi'_{v})\right) \right)  + O(s^2).
			\end{equation}
	 	
	 	\begin{definition}\label{def:Phi* local def}
	 		We define the modified  section $\Phi_{v}^*(g,s) \in I_{v}(g,s)$ by the formula	
				\begin{equation} 
	 				\Phi_{v}^*(g, s) \ := \ \Phi_{v}(g, s, \lambda_{\calV_{v}}(\varphi_{\calL,v})) - s \cdot  \left( \Phi \left(g, s, \lambda_{\calV_{v}}( \phi_{v})\right) +  \Phi\left(g, s, \lambda_{\calV'_{v}}(  \phi'_{v})\right)  \right).
	 			\end{equation}
	 	\end{definition}
		 
		Note that in general, 	the value at $s=0$ is given by
			\begin{equation} \label{eqn:Phi local value}
				\Phi_{v}^*(g, 0) =  \Phi_{v}(g, 0, \lambda_{\calV_{v}}(\varphi_{\calL,v})).
			\end{equation}
		Moreover,  for almost all $v$, we have that $	\Phi_{v}^*(g, s) $ is the standard spherical section, as per \Cref{rmk:modified local section unram}.
		
		For future use, we will gather some information about the Whittaker functionals $W_{m,v}(g, s, \Phi_{v}^*)$ near $s=0$.

		\begin{lemma} \label{lem:local Whittaker shift} 
			Suppose $\Phi(g,s) \in I_{v}(s, \chi)$ is an arbitrary section, $m \in F_{v}^{\times}$ and $t \in E_{v}^{\times}$. Then
				\begin{equation}
					W_{m,v}\left( m(t)^{-1}, s, \Phi \right) = |t|^{s-1}_{v} \chi(\overline t) W_{N(t)^{-1}m,v}(e, s, \Phi).
				\end{equation}
			where $N(t) = N_{E_{v}/F_{v}}(t)$. 
			\begin{proof}
				By definition, we have
				\begin{equation} \label{eqn:Whittaker t shift} 
				W_{m,v}\left( m(t)^{-1}, s, \Phi \right) = \int_{b \in F_{v}}\Phi\left(w n(b) m(t)^{-1}, s\right) \, \psi_p(-mb) \, db.
				\end{equation}
				A straightforward computation shows that 
				\begin{equation}
				w \, n(b) \, m(t^{-1}) \ = \ m({\overline t} ) \, w \, n(b N(t)).
				\end{equation}
				so
				\begin{equation}
					\Phi\left(w n(b) m(t)^{-1}, s\right) = \Phi\left(m(\bar t) \, w \, n(bN(t)), s\right) = |t|_{v}^{s+1} \chi(\bar t) \Phi\left( w \, n(b N(t)), s\right) .
				\end{equation}
				Substituting this equation into  \eqref{eqn:Whittaker t shift} yields
				\begin{equation}
					\begin{split}
						W_{m,v}\left( m(t)^{-1}, s, \Phi  \right) &= |t|_{v}^{s+1} \chi(\bar t) \int \Phi\left(w n(N(t)b)  , s\right) \, \psi_p(-mb) \, db \\
						& =  |t|_{v}^{s+1} \chi(\bar t) \int\Phi\left(w  n(b)  , s\right)  \psi_p(-m b/N(t)) \, 	\frac{db}{|t|^{-2}_{v}}   \\
						&= |t|_{v}^{s-1} \chi(\bar t) W_{N(t)^{-1}m, v}(e, s, \Phi)
					\end{split}
				\end{equation} 
				where in the second line we applied the change of variables $b \mapsto N(t)^{-1} b$. 
			\end{proof}
		\end{lemma}
		
		\begin{proposition} \label{prop:Phi* nonarch}
			Suppose $v$ is a non-split non-archimedean place of $F$, and let $m \in  F_{v}^{\times}$.
			\begin{enumerate}[(i)]
				\item $ W_{m,v}(e, 0, \Phi_{v}^*) = W_{m,v}(e, 0, \lambda_{\calV_{v}}(\varphi_{\calL,v}))$.
				\item If $V_{v}$ does not represent $m$, then
				\begin{equation}
				W_{m,v}(e, 0, \Phi^*_v) = 0 
				\end{equation}
				\item Suppose $V_{v}$ does not represent $m$.  Then
				\begin{equation}
					W'_{m,v}(e, 0, \Phi^*_{v})  - W'_{N(\varpi)^{-1}m, p}( e, 0, \Phi^*_{v}) 
					 = \ - \frac12 \log N(\lie p_E) \cdot W_{m,v}(e, 0, \lambda_{V'_{v}}(\varphi_{\calL}))
				\end{equation}
				where $\lie p_E$ is the prime ideal of $E$ inducing the valuation $v$ and $\varpi$ is a uniformizer of $E_v$. 
			\end{enumerate}
			\begin{proof}
				Part \textit{(i)} follows immediately from \eqref{eqn:Phi local value}. Part \textit{(ii)} follows from part \textit{(i)} and \Cref{prop:local whittaker dichotomoy}.
				
				For the third part, using \eqref{eqn:whittaker prime to -eps} and \eqref{eqn:rescale Schwartz fn} , we have  
					\begin{equation}
						W_{m, v}\left(e, 0, \lambda_{V'_{v}}(\varphi_L)\right) = |t|_{v} \chi(t) W_{m, v} \left(m(t)^{-1}, 0, \lambda_{V^{-\epsilon}_{v}	}(\varphi^0_{v}) \right)
				\end{equation}
				where $\Lambda = t \calO_{E,v}$ as in the argument leading to \eqref{eqn:whittaker prime to -eps}. Hence, by \Cref{lem:local Whittaker shift}, we conclude
					\begin{equation}  \label{eqn:Whittaker rel std 1}
						\begin{split}
									W_{m, v}\left(e, 0, \lambda_{V'_{v}}(\varphi_L)\right) &= \chi(t \overline t) W_{N(t)^{-1}m,v}\left( e, 0, \lambda_{V^{- \epsilon}_{v}}( \varphi^0_{v}) \right) \\
									& = W_{N(t)^{-1}m,v}\left( e, 0, \lambda_{V^{- \epsilon}_{v}}( \varphi^0_{v}) \right);
						\end{split}
					\end{equation}
				here we use the fact  $\chi(t \overline t)= 1$.  
				
				On the other hand, applying \Cref{lem:local Whittaker shift} again and the definition of $\Phi^*(g,s)$, we find
				\begin{equation} \label{eqn:Whittaker rel std 2}
					W_{m,v}(e, s,\Phi_{v}^*) = W_{N(t)^{-1} m, v}\left(e, s, \lambda_{V^{\epsilon}_{v} }(\varphi^0_{v}) \right) + O(s^2)
				\end{equation}
			   The Whittaker functions appearing on the right hand sides of \eqref{eqn:Whittaker rel std 1}  and \eqref{eqn:Whittaker rel std 2} were computed explicitly by Yang \cite[\S 2]{YangCM}, and the desired relations can be easily extracted from these formulas. 
			\end{proof}
		\end{proposition}
	
	The archimedean analogue is as follows. For later purposes, it will be convenient to introduce an additional parameter $y_v \in F_v = \bbR$ with $y_v >0$, which  will ultimately represent  the imaginary part of a variable in a Hilbert upper half-space on which our modular forms are defined.
	\begin{proposition} \label{prop:Phi* arch}
		Suppose $v | \infty$, and let 
			\begin{equation}
				\Phi^+_v(g,s) = \Phi(g,s, \lambda_{\calV_v}(\varphi_v)) \qquad {\text{and}} \qquad 				\Phi^{-}_v(g,s) =  \Phi(g,s, \lambda_{\calV'_v}(\varphi_v)).
			\end{equation}
		where $\varphi_v(x) = e^{-\pi \langle x, x \rangle_{\calV_v}}$.
		Recall here that $\calV_v$ and $\calV'_v$ have signatures $(1,0)$ and $(0,1)$, respectively. 
		
		For $m \in F^{\times}$ and $y_v \in F_v \simeq \bbR$ with $y_v > 0$, we set 
			\begin{equation}
				\calW_{m,v} (y_v, s,\Phi^{\pm}) := y_v^{-1/2} W_{m,v}( m(\sqrt{y_v}) , s,  \Phi^{\pm}_v ) \cdot e^{2 \pi y_v m_v}.
			\end{equation}
			where $m_v = \sigma_v(m)$ with $\sigma_v \colon F \to \bbR$ denoting the  real embedding corresponding to  the place $v$. 
			
			Suppose that $ m_v < 0$. Then for $t \in \bbR_{>0}$, we have  
				\begin{equation}
					\frac{d}{dt} \calW_{m,v}'(ty_v, 0, \Phi^+) \ = \ t^{-1} \calW_{m,v} (ty_v, 0,\Phi^-)
				\end{equation}
				and 
					\begin{equation}
						\lim_{t \to \infty} \calW_{m,v}(ty,0, \Phi^+) = 0.
					\end{equation}
	\begin{proof}
			This follows from the explicit formulas in \cite[Proposition 14.1]{KRYFaltings}.
	\end{proof}
	\end{proposition}
	
		At this stage, we have defined a local section $\Phi^*_{v}(s)$ for all non-split non-archimedian primes $v$. If $v | \infty$, we set $\Phi^*_v(g,s) = \Phi^+_v(g,s)$ as in the previous proposition. Finally, 	if $v$ is a split prime of $E$, we set $\Phi^*_{v}(g,s) = \Phi\left(g,s, \lambda_{V_{v}}(\varphi_L)\right)$
			
	Putting everything together, we make the following definition.
	\begin{definition} Let	
			\begin{equation}
				\Phi^*(s) := \otimes \Phi_v^*(s) \in I(s, \chi).
			\end{equation}
	\end{definition}

	\begin{remark}
		Let $\Sigma^{\bad} = \{ v_1, \dots, v_r\}$ denote the set of non-archimedean places  of $E$ such that either $E_v / F_v$ is ramified, or $E_v / F_v$ is an unramified field extension but $\calL_v$ is not self-dual. Recall that $\Phi^*_{v}(g,s) = \Phi(g, s, \lambda_{\calV_{v}}(\varphi_{\calL}))$ if $v \notin \Sigma^{\bad}$.
		
		For each $v_i \in \Sigma^{\bad}$, consider the collection $\{ \calV_w \ | \ w \neq v_i\} \cup {\calV'_{v_i}}$; this is a coherent collection, in the sense of \Cref{sec:prelim Hermitian}, so there exists an $E$-Hermitian space $W_i$ such that $W_{i,v_i} \simeq \calV'_{v_i}$ and $ W_{i,w} \simeq \calV_w$ at all other places. 
		
		A moment's contemplation of the definition of the local section $\Phi^*_{v}(g,s)$, cf.\ \Cref{def:Phi* local def}, shows that there exists Schwartz functions
			\begin{equation}
				\varphi_i \in S(W_i), \qquad \varphi_i' \in S(\calV)
			\end{equation}
		such that
			\begin{equation}
				\begin{split}
					\Phi^*(g,s) &= \Phi\left( g, s, \lambda_{\calV}\left( \varphi_{L} \right) \right)  + s \left[ \sum_{v_i \in \Sigma^{\bad}} \Phi(g, s, \lambda_{W_i}(\varphi_i)) +  \Phi(g, s, \lambda_{\calV}(\varphi_i')) \right] \\
					&\qquad \qquad  \qquad \qquad \qquad \qquad + O(s^2) 
				\end{split}
			\end{equation}
		Taking Eisenstein series, we find
			\begin{equation}
				E(g,s,\Phi^*) =  E \left( g, s, \lambda_{\calV}\left( \varphi_{L} \right) \right)  
				  + s \left[ \sum_{v_i  } E(g, s, \lambda_{W_i}(\varphi_i)) +  E(g, s, \lambda_{\calV}(\varphi_i')) \right]+ O(s^2) 
			\end{equation}
		The terms $E(g, s, \lambda_{\calV}(\varphi_i'))$ all vanish at $s=0$, since they are incoherent. Thus, at $s=0$ we have
			\begin{equation}
				E'(g, 0, \Phi^*) =  E' \left( g, 0, \lambda_{\calV}\left( \varphi_{L} \right) \right) + \sum_{v_i} E(g, 0, \lambda_{W_i}(\varphi_i)) .
			\end{equation}
		In other words $E'(g,0,\Phi^*)$ can be expressed as sum of the central derivative of the incoherent Eisenstein series $ E' \left( g, 0, \lambda_{\calV}\left( \varphi_{L} \right) \right)$, plus the central \emph{values} of  coherent Eisenstein series attached to  totally positive definite spaces.
		$\diamond$
	\end{remark}

	\subsection{Green functions and special cycles} \label{sec:Green fns}
		We now turn our attention to the construction of a family of elements $\widehat Z(m,\tau) \in \widehat \Pic_{\bbA}(\calM)$.	
		
		We begin with a local construction. Let $v$ be a non-archimedean place of $F$ that is non-split in $E$, and recall that we had fixed a local system $\bbV^{(v)} $ on $\calM_{/ E_{v}^{\alg}}$ of one-dimensional Hermitian spaces over $E$, as in \Cref{sec:data}.

		As part of our data, there is an   isometry 
		\begin{equation}
			\beta^{(v)}_z \colon	\bbV^{(v)}|_z \otimes \bbA_E \isomto \calV^{v} \times \calV'_{v}, 
		\end{equation}  
		at each point  $z \in \calM(E_{v}^{\alg})$. Pulling back the Schwartz function $ \varphi_{\calL, \bbA}$ along this map, we obtain a Schwartz function 
		\begin{equation} \label{eqn:varphi_z def}
			\varphi_z = \bigotimes_{w \leq \infty} \varphi_{z,w} \in S(\bbV^{(v)}|_z \otimes \bbA_E)
		\end{equation}
		on the adelification of the fibre of $\bbV^{(v)}$ at $z$.

		 For  $z \in \calM(E_{v}^{\alg})$ and $x \in  \bbV^{(v)}|_z $, we define a function
			\begin{equation}
				g(x, z) :=  \int_{\substack{ \lambda \in E_{v}^{\times} \\ |\lambda|_{v} \geq 1}}  \varphi_{z,v}(\lambda x) \, d^{\times} \lambda \cdot \log N_{E/{\bbQ}}(\lie p)
			\end{equation}
		where $d^{\times} \lambda$ is the multiplicative Haar measure on $E_{v}^{\times}$ such that the volume of $\calO_{E, v}^{\times}$ with respect to $d^{\times}\lambda$ is one, and $\lie p$ is the prime ideal of $E$ corresponding to the place $v$. 
		
		There are two main motivations for this construction. The first is that it strongly parallels Kudla's construction in the archimedean place, cf.\ \eqref{eqn:Green Arch} below. The more important point is that $g(x,z)$ naturally encodes the Arakelov degrees of divisors locally on $\calM$.
		
		We make the second point more precise. Suppose the point $z \in \calM(E_v^{\alg})$ is defined over a finite extension $L  $ of $E_{v}$,  and admits an integral extension $\widetilde z \in \calM(\calO_{L})$. Suppose further that there is a line bundle $\tilde \omega$ over $\calO_{L}$ and a map
			\begin{equation}
					s \colon \bbV^{(v)}|_z \otimes_{E_v} L \to \tilde \omega \otimes_{\calO_L} L;
			\end{equation}
		we may normalize $s$ to have the property that
			\begin{equation} \label{eqn: s_z normalization}
					s(x) \in \pi  \, \widetilde \omega \iff \varphi_{z, v}(x) = 1
			\end{equation}
		 for every $x \in \bbV^{(v)}|_z \otimes_E E_v$, where $\pi $ is a uniformizer of $L$. 
		 
		 For any $x$ with $s(x) \in \widetilde \omega$, its vanishing locus $Z(s(x))$ defines a divisor on $\Spec \calO_L$.	Let $I = (f) $ denote the (principal) ideal that cuts out $Z(s(x))$; in other words, $I$ is the largest ideal such that the image of $s(x)$ in $\widetilde \omega \otimes_{\calO_L} \calO_L/I$ vanishes. We define the normalized \emph{Arakelov degree} of the divisor $Z(s(x))$ to be 
		\begin{equation}
			\widehat \deg \, Z(s(x)) :=  \frac{1}{[L:E_{v}]}\mathrm{length}( \calO_{L} / I ) \, \log |\kappa(L)| = \frac{1}{[L:E_{v}]} \ord_{L}(f) \, \log |\kappa(L)|
		\end{equation}
		where $\kappa(L)$ is the residue field of $L$. Note that the Arakelov degree is unchanged if $L$ is replaced by a larger field. 
		
		\begin{proposition}\label{prop:Green Arakelov}
			With the notation as in the previous paragraph, we have
				\begin{equation}
					g(x, z) = \deghat Z(s(x))
				\end{equation}
			for all $x \in \bbV^{(v)}|_z$ with $s(x) \in \widetilde \omega$. 
			
			\begin{proof}
				If $\pi$ is a uniformizer of $L$, and $I = (f)$ is the ideal defining $Z(s(x))$, then $\ord_L(f)$ is, by construction, the largest integer such that $\pi^{- \ord_L(f)} s(x) \in \widetilde \omega$.  
				
				On the other hand, suppose $\varpi$ is a uniformizer of $E_{v}$ such that $\varpi= \pi^e$, with $e = e(L/ E_{v})$ the ramification index.  Then, using the normalization \eqref{eqn: s_z normalization}, we may write
				\begin{equation}
					\begin{split}
						g(x, z) &= \int_{\substack{\lambda \in E_{v}^{\times} \\ |\lambda|_{v} \geq 1}} \varphi_{z,v}(\lambda x)  \, d^{\times}\lambda  \cdot \log N(\lie p) \\
						& = \sum_{m \geq 0} \int_{\calO_{E,v}^{\times}}  \varphi_{z,v}(\varpi^{-m}u x)   d^{\times}u  \cdot \log  N(\lie p) \\
						&= \max \left( m \ | \ \varpi^{-m} s_x \in \widetilde \omega \right) \cdot \log  N(\lie p)  \\
						&= (\ord_{\pi}(f) / e) \cdot \log  N(\lie p) \\
						&= \frac{\ord_{\pi}(f)}{[L: E_{v}]} \cdot \log |\kappa(L)|
					\end{split}
				\end{equation} 	
				where in the last line, we used the identity $[L:E_{v}] = e(L / E_{v}) \cdot f(L/E_{v})$. This last expression is $ \widehat\deg \, Z(s(x))$, by definition.
			\end{proof}
		\end{proposition}
	In practice, one often has a global line bundle $\omega$, for which the vanishing locus of a section represents a moduli or deformation-theoretic problem; some examples of this phenomenon are presented in \Cref{sec:examples}.
	
	We record  another crucial property  of this construction: 
	\begin{lemma} \label{lem:Green equiv props}
		Suppose $z \in \calM(E_{v}^{\alg})$ and $x \in \bbV^{(v)}|_z$. For every $t \in H(\bbA_F)$, there is an isometry $r(t) \colon \bbV^{(v)}|_z \to \bbV^{(v)}|_{t \cdot z}$ 
			such that 	
			\begin{equation}
				\varphi_{z}(x) = \varphi_{t \cdot z}( t^{-1}r(t)  x)
			\end{equation}

		\begin{proof}
			This follows immediately from the equivariance property of $\beta^{(v)}$, as in \eqref{eqn:beta equivariance}.
		\end{proof}
	\end{lemma}
	
	Finally, for $\alpha \in F^{\times}$, we define a function $g_{v}(\alpha)$ on $\calM(E_{v}^{\alg})$ by the formula
		\begin{equation}
			\begin{split}
				g_{v}(\alpha) (z) &:= \sum_{\substack{x \in \bbV^{(v)}|_z \\ Q(x) = \alpha}} \varphi^{v}(x)  \cdot g_v(x,z) \\
				& =   \log N(\lie p)   \sum_{\substack{x \in \bbV^{(v)}|_z \\ Q(x) = \alpha}} \varphi^{v}_z(x) \int_{\substack{\lambda \in E_{v} \\ |\lambda|_{v} \geq 1 }} \varphi_{z,v}(\lambda x) \, d^\times \lambda
			\end{split}	
		\end{equation}
	where $\varphi^{v}_z(x) = \prod_{w \neq v} \varphi_{z,w}(x)$. 
	
	For later comparison to the Fourier expansion of a Hilbert modular form, it will be convenient to introduce an additional parameter $\tau \in \bbH_E$, where
		\begin{equation}	
			\mathbb H_E := \{  \tau = (\tau_v)_{v|\infty} \in \bbC^d \ | \ \mathrm{Im} (\tau_v) > 0 \text{ for all } v  \}
		\end{equation}
	is the Hermitian upper half space attached to $E$ (here $[E:\bbQ] = 2d$). For $\tau = x + iy \in \bbH_E$, choose a totally positive element $a \in F\otimes_{\bbQ} \bbR$ such that $a^2 = y$, and for an archimedean place $v$, set
			\begin{equation}
				g_{\tau,v}  = n(x_v) m(a_v) \in G(F_v)  
			\end{equation}
	and let
		\begin{equation} \label{eqn:gtau infinity}
			  g_{\tau,\infty}  = (g_{\tau, v_1}, \dots, g_{\tau, v_d})\in G(F_{\bbR}).
		\end{equation}
	We then  define the Green function
		\begin{equation} \label{eqn:padic Green final def}
			g_{v}(\alpha, \tau)(z) := N(y)^{-\frac12} \sum_{\substack{x \in \bbV^{(v)}|_{  z } \\ Q(x) = \alpha}}  \left( \omega_{\bbV^{(v)}_z}(g_{\tau,\infty}) \varphi_{z, \infty}(x) \right) \varphi_{z,f}^{v}(x)  \, g_v(x,z)
		\end{equation}
	where $\varphi_{z,\infty} = \prod_{w|\infty} \varphi_{z,w}$ and $N(y) =  \prod_{v|\infty} y_v   $.

	Now suppose $v$ is an archimedean place. For $z \in \calM(E^{\alg}_v)$ and $x \in \bbV^{(v)}_z$, we define
		\begin{equation} \label{eqn:Green Arch}
			g_v(x,z) = \int_1^{\infty} \varphi_{z,v}(tx) \, e^{2 \pi t^2 \sigma_v(Q(x))} \frac{dt}{t}
		\end{equation}
	where $Q(x) = \frac12 \langle x, x \rangle$ and $\sigma_v \colon F \to \bbR$ is the embedding induced by $v$. 	This definition coincides with the Green function introduced by Kudla in \cite{KudlaCentralDerivs}. 
	
	We incorporate the variable $\tau = (\tau_w)_{w|\infty} \in \bbH_E$ as follows. Let	
	\begin{equation}
		g_v(x, \tau_v)(z) := y_v^{\frac12} \, g_v \left( \sqrt{y_v} \, x  \right)(z) \, e^{ 2 \pi  \sigma_v(Q(x)) \tau_v }, 
	\end{equation}
	and for $\alpha \in F^{\times}$,  define a function $g_v(\alpha, \tau) \colon \calM(E_{v}^{\alg}) \to \bbR$ by the formula
		\begin{equation} 
				g_v(\alpha, \tau)(z) = N(y)^{-\frac12} \sum_{ \substack{x \in \bbV^{(v)}|_z \\ Q(x) = \alpha}}  g_v(x, \tau_v )(z)   \left( \prod_{\substack{w | \infty \\ w \neq v}} \omega_{\bbV^{(v)}_w}(g_{\tau,w})\varphi_{z,w}(x) \right)   \varphi_{z,f}(x).
		\end{equation}

	At this stage, we have constructed a function $g_v(\alpha, \tau)$ for every  non-split place  $v$.
	
	\begin{lemma}
	For a given $\alpha \in F^{\times}$, there is at most one place $v$ such that $g_{v}(\alpha,\tau)$ is not identically zero.
	\begin{proof}
		Suppose $g_{v}(\alpha, \tau)(z) \neq 0$ for some place $v$ and $z \in \calM(E_{v}^{\alg})$. Then there is a vector in  $\bbV^{(v)}|_z$ of norm $m$, so  $\bbV^{(v)}|_z$ is isometric to $(E, \alpha x \overline y)$. But $\bbV^{(v)}|_{z}$ and $\bbV^{(w)}|_{z'}$ are pairwise non-isometric for any $v \neq w$ and geometric points $z \in \calM(E_v^{\alg})$ and $z' \in \calM(E^{\alg}_w)$, hence there can be at most one such place $v$. 
	\end{proof}
	
\end{lemma}

	Finally, we arrive at the definition of the special cycles:
		\begin{definition}
				Let $\alpha \in F^{\times}$, and define the special cycle $\widehat Z(\alpha,y)$ in terms of the family of Green functions:
					\begin{equation}
						\widehat Z(\alpha,\tau) = \left( g_v(\alpha,\tau) \right)_{v \leq \infty} \in \widehat\Pic_{\bbA}(\calM_K)
					\end{equation}
		\end{definition}
	Of course, by the preceding lemma, there is at most one non-zero entry in the family.
	
\subsection{Proof of the main theorem }
	We are now ready to formulate our main theorem. 
	
	Let $E(g, s, \Phi^*)$ denote the Eisenstein series attached to the section $\Phi^*(g,s)$ as defined in \Cref{sec:Phi*}. It will be convenient to express $E(g,s, \Phi^*)$ in ``classical" coordinates, as follows: let $\tau \in \bbH_E$, and set
		\begin{equation}
			g_{\tau} = (g_{\tau, \infty}, 1, \dots ) \in G(\bbA_F)
		\end{equation} where $g_{\tau, \infty}$ is as in  \eqref{eqn:gtau infinity}. We then define
		\begin{equation}
			E(\tau,s, \Phi^*) := N(y)^{-\frac12} \, E(g_{\tau},s,\Phi^*),
		\end{equation}
	where
		\begin{equation}
			N(y) = \prod_{v|\infty} y_v;
		\end{equation}
		then $E(\tau,s,\Phi^*)$ is a (non-holomorphic) Hilbert modular form of parallel weight one. 
		
		Write the Fourier expansion of $E(\tau,s,\Phi^*)$ as
			\begin{equation}
				E(\tau, s, \Phi^*) = \sum_{\alpha \in F} E_{\alpha}(\tau, s, \Phi^*).
			\end{equation}
		
		\begin{theorem} \label{thm:main}
			For all $\alpha \in F^{\times}$, we have
			\begin{equation}
				\widehat \deg \, \widehat Z(\alpha,\tau) = -   \# \calM(\bbC)   \cdot E_{\alpha}'(\tau, 0, \Phi^*) .
			\end{equation}			
		
		\begin{proof} 
			By the product formula \eqref{eqn:Eis product}, we have  
				\begin{equation} \label{eqn:deriv Eis general}
					E_{\alpha}'(\tau, 0, \Phi^*) =N(y)^{-\frac12} \cdot \sum_{v \leq \infty} W'_{\alpha,v}(g_{\tau,v}, 0, \Phi^*_v) \prod_{w \neq v} W_{\alpha,w}(g_{\tau,w}, 0, \Phi^*_w).
				\end{equation}
			For $\alpha \in F^{\times}$, let 
				\begin{equation} \label{eqn:diff def}
					\mathrm{Diff}(\alpha) := \{v \leq \infty \mid \alpha \text{ is not represented by } \calV_v  \}.
				\end{equation}
			
			We observe that $\# \mathrm{Diff}(\alpha) \geq  1$. Indeed, if $\calV_v$ represents $\alpha$, then $\calV_v \simeq (E_v, \alpha xy')$. If $\mathrm{Diff}(\alpha) =\emptyset$, then this would hold true for all places, and so
				\begin{equation}
					\prod_{v \leq \infty} \inv_{v}(\calV_v) = 1,
				\end{equation}
			which   contradicts the assumption that $\calV$ is an incoherent collection.
			
			If $\# \mathrm{Diff}(\alpha) \geq 2$, then the vanishing criterion in \Cref{prop:local whittaker dichotomoy} implies that $E_{\alpha}'(\tau, 0, \Phi^*) = 0$. On the other hand, every space $\bbV^{(v)}|_{z}$ differs from $\calV$ at exactly one place, so no such space can represent $\alpha$; it follows that $\widehat Z(\alpha,\tau) = 0 $ in this case is well. 
			
			Thus, it remains to consider the case $ \mathrm{Diff}(\alpha) = \{ v \}$ is a singleton. Note that $v$ is necessarily non-split, so that $E_v = E \otimes_v F_v$ is a  field. 
			
			First, suppose $v$  is non-archimedean, and let $\lie p$ denote the prime ideal of $E$ inducing the valuation $v$ of $E$.  Then
				\begin{equation} \label{eqn:Eis prod nonarch}
				 N(y)^{\frac12} \,	E'_{\alpha}(\tau, 0, \Phi^*) =  W'_{\alpha,v}(e, 0, \Phi^*_{v}) \cdot W_{\alpha,\infty}(g_{\tau, \infty}, 0, \Phi^*_{\infty})\cdot \prod_{ \substack{w < \infty \\ w \neq v} } W_{\alpha,w}(e, 0, \Phi^*_w)  
				\end{equation}
			where 
				\begin{equation}
					W_{\alpha,\infty}(g_{\tau, \infty}, 0, \Phi^*_{\infty}) = \prod_{w|\infty} 	W_{\alpha,v}(g_{\tau, w}, 0, \Phi^*_{w}).
				\end{equation}
		On the other hand,  
				\begin{equation}
					\begin{split}
				 N(y)^{\frac12}  \, 	\widehat\deg  \,\widehat Z(\alpha,\tau) &=  N(y)^{\frac12}  \sum_{z \in \calM(E_{v}^{\alg})} g_{v}(\alpha,\tau)(z) \\
					&= \sum_{z \in \calM(E_{v}^{\alg})} \sum_{\substack{x \in \bbV^{(v)}|_z \\ Q(x) = \alpha}} \left( \omega_{\bbV^{(v)}_z}(g_{\tau, \infty}) \varphi_{z,\infty}(x ) \right)\varphi_{z,f}^{v}(x)  \\
					& \hspace{1.2 in} \times \int_{\substack{ \lambda \in E_{v}^{\times} \\ |\lambda|_{v} \geq 1}}  \varphi_{z,v}(\lambda x) \, d^{\times} \lambda  \, \log N(\lie p) 
					\end{split}
				\end{equation}
			with notation as in \eqref{eqn:padic Green final def}.  Next, for $\lambda \in E_{\lie p}^{\times}$, let 
				\begin{equation}
					\underline \lambda = (1, \dots, 1,  \underbrace{\lambda}_{\mathclap{\lie p \text{'th  component}}}, 1, \dots) \in \bbA_E^{\times}
				\end{equation}
			and write 
				\begin{equation}
					g = g_{\tau}  \, m(\underline \lambda) \in G(\bbA_F)
				\end{equation}
			so that 
				\begin{multline}
					\deghat \widehat Z(\alpha,\tau) = N(y)^{-\frac12 } \log N(\lie p) \\
					\times \int_{|\lambda|_{v} \geq 1}|\lambda|_{v}^{-1} \chi(\lambda)^{-1} \sum_{z \in \calM(E_{v}^{\alg})} \sum_{\substack{x \in \bbV^{(v)}|_z \\ Q(x) =\alpha}}  (\omega_{\bbV^{(v)}_z} (g) \varphi_z)(x)  \, d^{\times} \lambda.
				\end{multline}
			We compute the double sum in the integrand using the Siegel-Weil formula, as follows.

			Choose a  subgroup $K_H \subset H(\bbA_F)$  of finite index in $H(F) \backslash H(\bbA_F)$ such that $ \varphi_{\calL,\bbA}$ is $K_H$-stable, where $\varphi_{\calL, \bbA}$ is the Schwartz function \eqref{eqn:Schwartz L defn}. Assume furthermore that  the action of $H(F) \backslash H(\bbA_F)$ on $\calM(E_v^{\alg}) = \calM(E^{\alg})$ factors through the class group
				\begin{equation}
					\calC := H(F) \backslash H(\bbA_F) / K_H.
				\end{equation}
			
			Fix a set of representatives
				\begin{equation}
					z_1, \dots, z_r
				\end{equation}
			for the orbits of the action of $\calC$ on $\calM(E^{\alg})$, and for each $i$, let $c_i$ denote the order of the stabilizer of $z_i$ with respect to this action. 
			
			Then, for any $g \in G(\bbA_F)$, and $i = 1, \dots r$, we have
				\begin{equation}
					\begin{split}
					 \sum_{z \in \calM(E_{v}^{\alg})} \sum_{\substack{x \in \bbV^{(v)}|_z \\ Q(x) = \alpha}}  (\omega_{\bbV^{v}_z} (g) \varphi_z)(x) &= 
					 \sum_{i=1}^r \frac{1}{c_i}\sum_{t \in \calC}  \sum_{\substack{x \in \bbV^{(v)}|_{t\cdot z_i} \\ Q(x) =\alpha}} (\omega_{\bbV^{v}_{t \cdot z_i}} (g) \varphi_{t \cdot z_i})(x) \\
					 &=  \sum_{i=1}^r \frac{1}{c_i} \sum_{t \in \calC} \sum_{\substack{x \in \bbV^{(v)}|_{  z_i} \\ Q(x) = \alpha}} (\omega_{\bbV^{v}_{ z_i}} (g) \varphi_{  z_i})(t^{-1}x) \\
					 &= \sum_{i=1}^r \frac{1}{c_i} \sum_{t \in \calC} \Theta_{\alpha}(g, t, \varphi_{z_i})
					\end{split}
				\end{equation}
			where the second line follows from \Cref{lem:Green equiv props}, and $\Theta_{\alpha}(g, t, \varphi_{z_i})$ is the $\alpha$'th Fourier coefficient of the theta function defined in \eqref{eqn:theta function defn}. 
			
			Now by the Siegel-Weil formula (\Cref{prop:SW}),  the product formula \eqref{eqn:Eis product}, and  applying the definition of the Schwartz function $\varphi_{z_i}$, cf.\ \eqref{eqn:varphi_z def}, we can write
				\begin{equation}
					\begin{split}
					\sum_{t \in \calC} \Theta_{\alpha}(g, t, \varphi_{z_i}) &= 
					\frac{1}{\mathrm{vol}(K_H, dh) } \int_{[H]} \Theta_{\alpha}(g, h, \varphi_{z_i}) dh \\
					&= \frac{1}{2 \mathrm{vol}(K_H)} E_{\alpha}(g, 0, \lambda_{\bbV^{v}_{z_i}}(\varphi_{z_i})) \\
					&= \frac{1}{2 \mathrm{vol}(K_H)}   \cdot \prod_{\substack{w \leq \infty \\ w \neq v}} W_{\alpha,w}(g_w, 0, \lambda_{\calV_w}(\varphi_{\calL,w})) \cdot  W_{\alpha,v}(g_{v}, 0, \lambda_{\calV'_{v}}\left(\varphi_{\calL, v} ) \right).
					\end{split}
				\end{equation}
		Note that this is independent of $i$. Moreover, we observe that
			\begin{equation}
				\# \calM(E^{\alg}) =  \# \calC \sum_{i=1}^r c_i^{-1}   =\left( \frac{1}{\mathrm{vol}( K_H) } \int_{[H]} dh \right)   \sum_{i=1}^r c_i^{-1}   = \frac{1}{\mathrm{vol}( K_H)} \sum_{i=1}^r c_i^{-1};
			\end{equation}
		thus, setting 
			\begin{equation}
				C := \frac{\# \calM(E^{\alg})}{2}
			\end{equation}
		and taking $g = g_{\tau, \infty} \, m(\underline \lambda)$, we have  
			\begin{multline} \label{eqn:deg product non-arch}
				\widehat \deg \, \widehat Z(\alpha,\tau) =C \cdot N(y)^{-\frac12} \cdot W_{\alpha, \infty}(g_{\tau, \infty}, 0,  \lambda_{\calV_{\infty}}(\varphi_{\calL, \infty}))\prod_{\substack{w < \infty \\ w \neq v}} W_{\alpha,w}(g_w, 0, \lambda_{\calV_w}(\varphi_{\calL,w})) 
				  \\ \times \log N(\lie p) \int_{\substack{ \lambda \in E_{v}^{\times} \\ |\lambda|_{v} \geq 1}} |\lambda|_{v}^{-1} \chi(\lambda)^{-1} W_{\alpha, v}(m(\lambda), 0, \lambda_{\calV'_{v}}(\varphi_{\calL, v})) \, d^{\times} \lambda.
			\end{multline}
		Now applying \Cref{lem:local Whittaker shift}, we may rewrite the second line above as
			\begin{equation} \label{eqn:whittaker mellin integral}
				\begin{split}
					\log N(\lie p) \int_{\substack{ \lambda \in E_{v}^{\times} \\ |\lambda|_{v} \geq 1}} |\lambda|_{v}^{-1} \chi(\lambda)^{-1} & W_{\alpha, v}(m(\lambda), 0, \lambda_{\calV'_{v}}(\varphi_{\calL, v})) \, d^{\times} \lambda \\
					 &=  	\log N(\lie p)  \int_{\substack{ \lambda \in E_{v}^{\times} \\ |\lambda|_{v} \geq 1}} W_{N(\lambda)\alpha, v}(e, 0, \lambda_{\calV'_v}(\varphi_{\calL, v})) \, d^{\times} \lambda.
				\end{split}
			\end{equation} 
		Note that since $\varphi_{\calL, v}$ is the characteristic function of a lattice,  
			\begin{equation}
				W_{u\alpha, v}(e, 0, \Phi_{v}) = W_{\alpha , v}(e, 0,\lambda_{\calV'_v}(\varphi_{\calL, v}))
			\end{equation}
		for all $u \in \calO_{F, v}^{\times}$. Fixing a uniformizer $\varpi$ of $E_{\lie p}$ and applying \Cref{prop:Phi* nonarch}(iii), we have
			\begin{equation}
				\begin{split}
					\eqref{eqn:whittaker mellin integral} &= \log N(\lie p)   \sum_{k \geq 0} 	W_{N(\varpi)^{-k}\alpha, v} \left( e, 0, \lambda_{\calV_{v}'}(\varphi_{\calL, v}) \right) \\
					&= -2 \, W_{\alpha,v}'(e, 0, \Phi^*_{v}).
				\end{split}
			\end{equation}
		As for the remaining places, \Cref{prop:Phi* nonarch}(i) gives
			\begin{equation}
				W_{\alpha,w}(e, 0, \lambda_{\calV_w}(\varphi_{\calL,w})) = W_{\alpha,w}(e, 0, \Phi^*_w) 
			\end{equation}
		for non-archimedean places $w\neq v$ , and, following definitions, we have
			\begin{equation}
				W_{\alpha,w}(g_{\tau, w}, 0, \lambda_{\calV_w}(\varphi_{\calL,w})) = W_{\alpha,w}(g_{\tau,w}, 0, \Phi^*_w)
			\end{equation}
		for $w| \infty$. 
		
		Substituting these expressions into \eqref{eqn:deg product non-arch} and comparing with \eqref{eqn:Eis prod nonarch} gives the result in this case.
		
		The case that $\mathrm{Diff}(m) = \{ v \}$ for an archimedean place $v$ is entirely analogous, and is left to the reader. 
		\end{proof}
	
		\end{theorem}
	
	\begin{remark} The theorem above readily generalizes to the situation where the compact open subgroup $K \subset T(\bbA_{\bbQ,f})$ is not neat.  Indeed, fix a neat subgroup $K_0 \subset K$, let $\calM_0$ denote the corresponding zero-dimensional Shimura variety, and let $\calM = [ \Gamma \backslash \calM_0]$ denote the stack quotient, where $\Gamma = K / K_0$. Let $j \colon \calM_0 \to \calM$ denote the natural map.
		
	Suppose $\widehat \omega$ is an adelic line bundle on $\calM$, given by a family of Green functions $(g_v)_v$ as in \Cref{sec:adelic metric}. If we define its arithmetic degree on $\calM$ to be
		\begin{equation}
					\deghat\, \widehat \omega = \sum_{v \leq \infty} \sum_{z \in \calM(E_v^{\alg})} \frac{1}{\# \Aut(z)}g_v(z),
		\end{equation}
	then one may verify that
		\begin{equation}
			 \frac{	\deghat  j^* \widehat \omega }{\# \calM_0(\bbC)} = \frac{\deghat \widehat\omega}{\# \calM(\bbC)}.
		\end{equation}
	where the left hand side is the ``stacky" cardinality. 
	This in turn implies that the main theorem continues to hold in this case. 
	
	Alternatively, one can take the definition of the arithmetic degree above, and make the minor changes required in the proof of \Cref{thm:main} directly.  \hfill $\diamond$
	\end{remark}
	
		\section{Examples} \label{sec:examples}
			\subsection{The weight  one Eisenstein series of Kudla-Rapoport-Yang}

				Suppose that $E \subset \bbC$ is an imaginary quadratic field, and let $\calO_E \subset E$ denote the ring of integers.
				 Let $\calM$ denote the moduli stack over $\Spec (\calO_E)$ whose $S$-points, for a connected $\calO_E$-scheme $S$, parametrizes pairs $(A, i)$ where $A$ is an elliptic curve over $S$, and $i \colon \calO_E  \to \End(A)$ is an action of $\calO_E $ such that the induced action on $\Lie(A)$ coincides with  the structural morphism $\calO_E \to \calO_S$  . Then $\calM$ is a zero-dimensional Shimura variety with
				 	\begin{equation}
					 		\calM(\bbC) \simeq \left[ E^{\times} \big\backslash \bbA_{E,f}^{\times} / \widehat{\calO_E}{}^{\times} \right].
				 	\end{equation}
				
				 In \cite{KRY-IMRN}, the authors assume that $E = \bbQ(\sqrt{-p})$ for a prime $p$,  construct a family of special cycles on $\calM$ via a moduli problem, and prove that their arithmetic degrees coincide with the Fourier coefficients of an incoherent weight one Eisenstein series. In this section, we explain how their results can be viewed in the context of the present work. 
				 
				 To begin, we fix our incoherent collection $\calV = (\calV_v)_v$ as follows. If $v$ corresponds to a rational prime $p$, we set $\calV_v = (E \otimes_{\bbQ} \bbQ_p, - 2x \overline y)$, while at the infinite place, we set $\calV_{\infty} = (\bbC, 2x \overline y)$; here we view $\bbC \simeq E \otimes_{\bbQ}\bbR$ via the fixed embedding $E \to \bbC$. 
				 As for the lattice $\calL$, we set $\calL_p = \calO \otimes_{\bbZ} \bbZ_p$ for each prime $p$. 
				 
				 Next, let $p$ denote a rational prime that is non-split in $E$, and take a geometric point $z \in \calM(E_p^{\alg})$  corresponding to an elliptic curve $A_z$ with its $\calO_E$-action. Let $A_{\overline z} \in \calM(\bbF_p^{\alg})$ denote the reduction mod $p$, which is equipped with an induced $\calO_E$-action $\iota\colon \calO_E \to \End(A_{\overline z})$, and set
				 	\begin{equation}
					 	\bbV^{(p)}|_z  := \{  x \in \End^0(A_{\overline z})  \mid \iota(a) \circ x = x \circ \iota (\bar a) \ \text{for all } a \in \calO  \};
				 	\end{equation}
				 this is the space of ``special quasi-endomorphisms" of $A_{\overline z}$ as in \cite{KRY-IMRN}. 
				 
				 Since $p$ does not split in $E$, the elliptic curve  $A_{\overline z}$ is necessarily supersingular, and so $\End^0(A_{\overline z})$ is a rational quaternion algebra that is non-split precisely at $\infty$ and $p$. Since $\End^0(A_{\overline z})$ contains an embedded copy of $E$, there is a decomposition
				 	\begin{equation}
				 		\End^0(A_{\overline z}) = \iota(E) \oplus \iota(E) \theta
				 	\end{equation}
				where $\theta$ is an element satisfying $\theta \iota (a) = \iota(\overline a) \theta$ for all $a \in \calO_E$. In particular, $\bbV^{(p)}|_z = \iota(E) \theta$ is a one-dimensional vector space over $E$. We equip $\bbV^{(p)}|_z$ with the E-Hermitian form $\langle \cdot , \cdot \rangle_z $ such that $\frac12 \langle x, x \rangle_z = \mathrm{Nrd}(x) = - x^2 $; here $\mathrm{Nrd}(x)$ is the reduced norm. 
				
				 We need to show that as $z$ varies, these spaces form an $H(\bbA_{\bbQ})$-equivariant local system, where $H(\bbA_{\bbQ}) = \{ t \in \bbA_E^{\times} \mid N(t) = 1 \}$. This is best seen via Serre's tensor construction as follows. For $t \in H(\bbA_{\bbQ})$, let $\lie a = t \widehat \calO_E \cap E$ denote the corresponding fractional ideal. Then we may identify $A_{t\cdot z} \simeq A_z \otimes_{\calO_E} \lie a$, and we obtain a natural $E$-linear quasi-isogeny
				 	\begin{equation}
					 	\varphi_t \colon A_{z} \to A_{z} \otimes_{\calO_E} \lie a \simeq A_{t\cdot z}, \qquad a \mapsto a \otimes 1.
				 	\end{equation}
				 This induces the desired isometry
				 	\begin{equation}
					 	r(t) \colon \bbV^{(p)}|_{t \cdot z} \to \bbV^{(p)}|_z, \qquad x \mapsto \varphi_t \circ x \circ \varphi_t^{-1}
				 	\end{equation}
				 of $E$-Hermitian spaces, where, abusing notation, we use the same symbol $\varphi_t$ to denote the induced quasi-isogeny between the special fibres of $A_z$ and $A_{t\cdot z}$.
				 
				 Now let $z \in \calM(E_p^{\alg}) $ as above, and let
				 	\begin{equation}
					 	\bbL^{(p)}|_z := \bbV^{(p)}|_z \cap \End(A_{\overline z}).
				 	\end{equation}
				 	At a prime $q \neq p$, we may identify the Tate module $\Ta_q(A_{\overline z}) \simeq \calO_{E,q}$, and so we obtain  isomorphisms
				 	\begin{equation}
					 	\End(A_{\overline z}) \otimes_{\bbZ} \bbZ_q \simeq \End(\Ta_q(A_{\overline z})) \simeq \End(\calO_{E,q}).
				 	\end{equation}
				The space of $\calO_E$-antilinear endomorphisms of $\calO_{E,q}$ is generated, as an $\calO_{E,q}$-module, by the  endomorphism $a \mapsto \bar a$; pulling this back through the above identifications, we obtain a generator $f \in \bbL^{(p)}|_z \otimes \bbZ_q$ with $f^2 = - 1$. This choice of basis vector induces an isometry $\bbL^{(p)}|_z \otimes \bbZ_q \simeq \calL_q$, or put differently,  an isometry 
					\begin{equation}
						\beta^{(p)}_{z,q} \colon \bbV^{(p)}|_{z} \otimes_{\bbQ} \bbQ_q \isomto \calV_q
					\end{equation} 
				identifying $\bbL^{(p)}_{z} \otimes \bbZ_q$ with $\calL_q$. 
				
				At the prime $p$,  note that we may choose the element $\theta $ such that its image in $\bbL^{(p)}|_z \otimes \bbZ_p$ is a generator, so that 
					\begin{equation}
						\bbL^{(p)}|_{z} \otimes \bbZ_p \simeq (\calO_{E,p}, 2 \mathrm{Nrd}(\theta) \, x \overline y)
					\end{equation}
				On the other hand, because $\End^0(A_{\overline z}) \otimes \bbQ_p$ is a division algebra, $\bbV^{(p)}|_{z,p} $ cannot represent $-2$, so $2 \mathrm{Nrd}(\theta) \equiv -2\varsigma_p \in \bbQ_p^{\times} / N(E_p^{\times})$, where $\varsigma_p$ is the element chosen in \Cref{sec:lin alg}. Moreover, $\End(A_{\overline z}) \otimes \bbZ_p$ is the unique maximal order, consisting precisely of  those endomorphisms whose reduced norms have non-negative $p$-adic valuation. Hence $\ord_p(\mathrm{Nrd}(\theta)) = 0$ or $1$, depending on whether $p$ is ramified or inert, respectively, and we find $\mathrm{Nrd}(\theta) = - N(u) \varsigma_p$ for some $u \in \calO_{E,p}^{\times}$. To summarize, we have an isometry
					\begin{equation}
							\beta^{(p)}_{z,q} \colon \bbV^{(p)}|_{z} \otimes \bbQ_p \isomto \calV'_q
					\end{equation}
				identifying $\bbL^{(p)}|_z \otimes_{\bbZ} \bbZ_p $ with $\calL_p$. 
				Finally, at the archimedean place, we have that both $\bbV^{(p)}_{z, \infty}$ and $\calV_{\infty}$ are positive definite one dimensional Hermitian spaces, and so we may fix an isometry $\beta^{(p)}_{z, \infty}$ between them.
				
				Putting these local factors together, we obtain an isometry
					\begin{equation}
						\beta^{(p)}|_z \colon \bbV^{(p)} \otimes_E \bbA_E \isomto \calV^p \times \calV'_p
					\end{equation}
				which at the non-archimedean places identifies $\bbL^{(p)}|_z \otimes \widehat \calO_E$ with $\widehat \calL$; note that this latter condition determines the finite part of $\beta^{(p)}_z$ up to multiplication by an element of $\widehat \calO_E^{\times}$ of norm one. One can then  check that for $t \in H(\bbA_{\bbQ})$, the corresponding maps $\beta^{(p)}_z$ and $\beta^{(p)}_{t \cdot z}$ satisfy the required equivariance property \eqref{eqn:beta equivariance}. 
				
				Now consider the archimedean place $v = \infty$. For a point $z \in \calM(E_v^{\alg}) = \calM(\bbC)$, the homology  group $H_1(A_z(\bbC), \mathbb Z)$ has an induced $\calO_E$ action, and we define
					\begin{equation}
						\bbL^{(\infty)}|_z = \{  x \in \End \left( H_1(A_z(\bbC), \mathbb Z)  \right) \ | \ x \circ \iota(a) = \iota(\overline a) \circ x \} 
					\end{equation}
				and
					\begin{equation}
						\bbV^{(\infty)}|_z  = \bbL^{(\infty)}|_z \otimes_{\bbZ} \bbQ,
					\end{equation}
				equipped with the $E$-hermitian form determined by $ \frac 12 \langle x,x\rangle_z = - x^2$. A similar argument to the non-archimedean case reveals that there is an isometry
					\begin{equation}
						\beta^{(\infty)} \colon \bbV^{(\infty)}|_z \otimes_{E} \bbA_E \isomto \calV^{\infty} \times \calV_{\infty}'
					\end{equation}
				satisfying the equivariance \eqref{eqn:beta equivariance}.
				
				At this point, we have specified the data $(\calV, \widehat \calL, (\bbV^{(v)}, \beta^{(v)} ))$ as set out in \Cref{sec:data}.
				Carrying out the construction in \Cref{sec:Green fns}, we obtain a family of Green functions $g_v(\alpha, \tau) \colon \calM(E_v^{\alg}) \to \bbR$, and consequently a family of special cycles $\widehat Z(\alpha, \tau) \in \widehat\Pic_{\bbA}(\calM)$ for $\alpha \in \bbQ^{\times}$ and $\tau \in \bbH$. 
				
				On the other hand, in \cite{KRY-IMRN}, the authors define a family of special cycles $\mathscr Z(\alpha)$ via a moduli problem. When $\alpha>0$, the link between the two constructions is given by \Cref{prop:Green Arakelov}. To make this more precise, recall that $\mathscr Z(\alpha)$ is defined by the following moduli problem: for a scheme $S$ over $\Spec(\calO_E)$, the $S$-points parametrize tuples $(A, i, x)$ where $(A,i)$ is a point of $\calM(S)$ and $x \in \End_S(A)$ is an endomorphism satisfying $x \circ i(a) = i(\overline a) \circ x$ for all $a \in \calO_E$ and $- x^2 =  \alpha$. 
				
				Now suppose we have a point $\overline{z} \in \calM(\bbF_p^{\alg})$, corresponding to a pair $(A_{\overline z}, i_{\overline z})$.  By Gross' theory of canonical liftings, this pair admits a canonical lift to characteristic zero; more precisely, choosing a sufficiently large finite extension $L$ of $\bbQ_p$, we have a (canonical) lift to a pair $(A_{\tilde z}, i_{\tilde z})$ over $\calO_L$.   
 				The algebraic de Rham homology group $H_1^{dR}(A_{\tilde z})$ is a free $\calO_L$ module of rank 2, equipped with a Hodge filtration
						\begin{equation}
							\begin{CD} 
							0 @>>> \Fil^1(A_{\tilde z}) @>>> H_1^{dR}(A_{\tilde z})  @>>> \Lie(A_{\tilde z}) @>>> 0.
							\end{CD}
						\end{equation}
				Furthermore, we may canonically identify $H_1^{dR}(A_{\tilde z}) $ as the value of the Grothendieck-Messing crystal of $A_{\overline z}$ at $\calO_L$; in particular, any element $x \in \bbL^{(p)}|_{z}$ induces an endomorphism $\tilde x$ of $H_1^{dR}(A_{\tilde z})$. 
				Consider the line bundle
						\begin{equation}
							\widetilde \omega = \Hom(\Fil^1(A_{\tilde z}), \Lie(A_{\tilde z})) .
						\end{equation}
			  Given $x \in \bbL^{(p)}|_{z}$, we obtain a section $s(x) \in \widetilde  \omega  $ by assigning to $x$ the composition
			  	\begin{equation}
				  s(x) \colon	\Fil^1(A_{\tilde  z})   \stackrel{\tilde x}{\longrightarrow} H^{dR}_1(A_{ \tilde z})   \to \Lie(A_{\tilde z}) 
			  	\end{equation}
		  	 where the second arrow is the natural projection.
			 By Grothendieck-Messing theory, the deformation locus of $x$ is precisely the vanishing locus of $s(x)$. In particular, since $x \in \bbL^{(p)}|_z$ is an endomorphism of $A_{\overline z}$,  we have $s(x) \in \pi \widetilde \omega$, where $\pi $ is a uniformizer of $\calO_L$. Tensoring with $\bbQ$, we obtain a map $s = s_z \colon \bbV^{(p)}|_z \to \tilde \omega \otimes \bbQ$ satisfying the compatibility condition \eqref{eqn: s_z normalization}. 
			
			Now suppose $\overline y \in \mathscr Z(\alpha)(\bbF_p^{\alg})$ lies above the point $\overline z \in \calM(\bbF_p^{\alg})$ so that $\overline y$ corresponding to a tuple $(A_{\overline z}, i_{\overline z}, x)$ with $x^2  = - \alpha$. In light of the above observations, \Cref{prop:Green Arakelov} implies that
				\begin{equation}
					 g_v(x, z) =  \deghat Z(s(x)) = \log  \# (\calO_{\mathscr Z(\alpha), \overline y})
				\end{equation}
			and it follows that for $\alpha >0$, we have
				\begin{equation}
					\deghat \widehat Z(\alpha, \tau) = \deghat \mathscr Z(\alpha) \, q^{\tau}.
				\end{equation} 
				
			When $\alpha<0$, a direct computation shows that the above identity continues to hold, where $\mathscr Z(\alpha) = \mathscr Z(\alpha,\tau)$ is defined in \cite[\S 6]{KRY-IMRN}.
							 
			Thus,  our main theorem \Cref{thm:main} specializes to the following identity, which is due to Kudla, Rapoport and Yang  in the case that the discriminant of $E$ is prime, and generalized to arbitrary discriminant by Kudla and Yang. 			
			\begin{theorem}[{\cite{KRY-IMRN,KudlaYang}}] 
					For any $\alpha \neq 0$, we have 
						\[
								\deghat \mathscr Z(\alpha) \, q^{\alpha} = -  \# \calM(\bbC) \cdot E'_{\alpha}(\tau, 0, \Phi^*).
						\]
					\qed
			\end{theorem}

			\subsection{The Bruinier-Kudla-Yang formula}
				In this section, we revisit some aspects of a paper by Andreatta, Goren, Howard, and Madapusi Pera \cite{AGHMP} in which they prove an averaged version of Colmez' conjecture; at the heart of their proof is a formula (originally conjectured by Bruiner, Kudla, and Yang)   for the arithmetic intersection number  of a special divisor and a ``big CM cycle" on an orthogonal Shimura variety. Our aim is to highlight the relevant parts of their setup in the context of our discussion.

				Suppose that  $E$ is a CM field with $[E:\bbQ] = 2d$ and let  $F \subset E$ the maximal totally real subfield. 
				We define a zero-dimensional Shimura variety as follows; note that our notation here will differ slightly from previous sections, in order to align more closely with \emph{loc.\ cit}. Consider the tori
					\begin{equation}
						T_E  := \Res_{E/\bbQ} (\bbG_m), \qquad  T_F   := \Res_{F/\bbQ}(\bbG_m)  , 
					\end{equation}
				and let 	
					\begin{equation}
						T := T_E / \ker( N_{F/\bbQ} \colon T_F \to \bbG_m) .
					\end{equation}
				We have $T_E(\bbC) \simeq \prod_{\sigma \colon E \to \bbC } (\bbC_{\sigma})^{\times}$, where $\bbC_{\sigma} := E \otimes_{E, \sigma} \bbC$. Choose an embedding $\sigma_1 \colon E \to \bbC$, and let $h_E \colon \mathbb S \to T_{E,\bbC}$ denote the homomorphism given on $\bbC$-points by the formula $h_E(z_1,z_2) = (z_1,  z_2,1, \dots, 1)$ where $z_1 \in \bbC_{\sigma_1}$ and $z_2 \in \bbC_{\overline \sigma_1}$, and all other components are one; we then define $h \colon \mathbb S \to T_{\bbR}$ to be the composition of $h_E$ with the projection to $T_{\bbR}$. 
				The reflex field of the pair $(T, h)$ is the field $\sigma_1(E) \subset \bbC$.

				Given a compact open subgroup $K \subset T(\bbA_f)$, the Shimura datum $(T,h)$ described above gives rise to   a zero-dimensional Shimura variety $Y_K$ defined over $E$, with 
					\begin{equation}
						Y_K(\bbC_{\sigma_1}) = T(\bbQ) \backslash \{ h \} \times T(\bbA_f) / K.
					\end{equation}
				
				Our next step is to specify the requisite data, as in \Cref{sec:data}. Let $\tau_1, \dots, \tau_d \colon F \to \bbR$ denote the set of real embeddings. 
				Fix an element $\xi \in F^{\times}$ such that $\tau_1(\xi) < 0$ and $\tau_i(\xi) > 0$ for $i = 2, \dots, d$.
				
				We define our incoherent family $\calV = (\calV_v)_{v}$ as follows: if $v$ corresponds to the embedding $\tau_1$, we set $\calV_v = (E_v, - \xi x \overline y)$, while at every other place $v$, we set $\calV_v = (E_v, \xi x \overline y)$. 
				In addition, we fix an $\widehat\calO_E$-stable, $K$-invariant lattice $\widehat \calL \subset \prod_{v < \infty} \calV_v$.
				
				Next, we need to construct the family of local systems $\{  ( \bbV^{(v)}, \beta^{(v)})\}$, which will require a little bit of preparation.  Let $(V, Q)$ denote the rational quadratic space of signature $(2d-2,2)$ whose underlying space is $V = E$ with quadratic form $Q(x) = tr_{F/\bbQ}(\xi x \overline x)$. As explained in \cite[\S 5]{AGHMP}, this data gives rise to an embedding
				\begin{equation}
				j_{/ \bbQ} \colon Y_K \to \calM
				\end{equation}
				where $\calM$ is (the base change to $\Spec(\calO_E)$ of) an integral model of the  $\GSpin$ Shimura variety associated to $V$; note that this model depends on the choice of a maximal lattice $\widehat{L}^{\max} \subset V$ containing $\mathscr L$, which we suppress from the notation. This extends to a morphism 
					\begin{equation}
						j \colon \calY_K \to \calM
					\end{equation}
			where $\calY_K$ is the normalization of $Y_K$ in $\calM$. 
				
				The variety $\calM$ comes with an abelian scheme $\calA^{KS} \to \calM$,  called the Kuga-Satake abelian scheme, which is of relative dimension $2^{d-1}$. Pulling back along $j$, we obtain an abelian scheme
				\begin{equation}
				\calA := j^* \calA^{KS}  
				\end{equation}
				over $\calY_K$. For a connected $\calY_K$-scheme $S$, let
				\begin{equation}
				L(\calA_S) 	\subset \End(\calA_S)
				\end{equation}
				denote the space of \emph{special endomorphisms}, as defined in \cite[\S 4.3, \S4.5]{AGHMP}  The precise definition is quite technical, so we content ourselves here with a summary of the relevant properties, beginning with the fact that $L(\calA_S)$ admits the structure of an $\calO_E$-module, together with a positive definite Hermitian form. 
				
				It will be convenient to work with ``good" places, as in \cite[Definition 5.3.3]{AGHMP}. Let $v$ denote a non-split, non-archimedean place of $F$. We say that $v$ is good if the level structure $K$ is maximal at $v$ and either \textit{(a)} the extension $E_v/F_v$ is  unramified,  the lattice $\calL_{v}$ is self-dual, and $ L^{\max}_v = \calL_{v} \subset E_v$; or  \textit{(b)} if $E_v/F_v$ is ramified and  $\calL_ {v}  \subset  L^{\max}_v \subsetneq \lie d_{E_{v}/F_{v}}^{-1} \calL_{v}$. Note all but finitely many places are good. 
				
				Now suppose $v$ is a good (non-archimedean) place, and let $y \in \calY_K(E_v^{\alg})$ denote a geometric point; by construction, there is unique lift $\tilde y \in \calY_K(\calO_{v}^{\alg})$, and a reduction $\overline y \in \calY_K(\kappa(v)^{\alg})$, where $\kappa(v)^{\alg}$ is the residue field of $E_v$. We define the local system $\bbV^{(v)}$ by setting
				\begin{equation}
				\bbV^{(v)}|_y :=    L(\calA_{\overline y}) \otimes_{\bbZ} \bbQ. 
				\end{equation}
				By \cite[Corollary 5.4.6, Proposition 7.6.2]{AGHMP}, this space  is one-dimensional over $E$, and  there is an isometry 
				\begin{equation}
				\beta^{(v)}_y \colon \bbV^{(v)}|_y \otimes_E \bbA_E \isomto \calV^{v} \times \calV'_v
				\end{equation} 
				that identifies $L(\calA_{\overline y})_w$ with   $L^{\max}_w$ at all non-archimedean places $w$. In addition, \cite[Lemma 7.6.3]{AGHMP} implies that the  map  $\beta^{(v)} \colon \bbV^{(v)} \otimes \bbA_E \to \underline{ \calV^{v} \times \calV'_v}$ of local systems is $H(\bbA_F)$-equivariant in the sense of \cref{eqn:beta equivariance}. 
				
				Applying the construction in \Cref{sec:Green fns}, we obtain a Green function $g_v(\alpha,\tau)$ for each good non-archimedean place $v$, and consequently a special cycle $\widehat Z(\alpha, \tau) \in \widehat \Pic_{\bbA}(Y_K)$ for $\alpha \in F^{\times}$ for any $\alpha$ with $\mathrm{Diff}(\alpha) = \{ v\}$ for some good place $v$,  cf.\ \cref{eqn:diff def}. 
				
				As with the previous example, the cycles in \cite{AGHMP} are defined in terms of a moduli problem, and are related to the cycles $\widehat Z(\alpha, \tau)$ via \Cref{prop:Green Arakelov}. To make this connection precise, for $\alpha \in F^{\times}$ ,  consider the moduli problem $\mathscr Z(\alpha)$ over $\calY_K$ whose points, for a connected $\calY_K$-scheme $S$,  are given by
				\begin{equation}
				\mathscr  Z(\alpha)(S)  = \left\{ x \in L(\calA_S) \ | \ \langle x, x \rangle = \alpha   \right\}.
				\end{equation}
				One has that  $\mathscr Z (\alpha) = 0$ if $\alpha$ is not totally positive or $\# \mathrm{Diff}(\alpha) \geq 2$, and that when $\mathrm{Diff}(\alpha) = \{ v\}$, the  moduli problem $\mathscr Z(\alpha) $ defines a zero-cycle that is contained in the special fibre of $\calY_K$ at $v$, 				see \cite[\S 7.6]{AGHMP} for details.
				
				To relate this moduli problem to the Green functions constructed in the present paper, we first observe that  the space $V$  (viewed as a representation $T$) together with the choice of lattice $L^{\max}$, gives rise to a filtered vector bundle $(V_{dR}, \Fil^{\bullet}(V_{dR}))$ over $\calY_K$, and moreover
					\begin{equation}
						\omega := \Fil^1(V_{dR})^{\vee}
					\end{equation}
				is  a line bundle; see \cite[Proposition 3.8.1]{AGHMP}. Now suppose  $\mathrm{Diff}(\alpha) = \{ v\}$ is a good place and let $\tilde y \in \calY_K(\calO_{E_v^{\alg}})$ with special fibre $\overline y \in \calY_K(\kappa(v)^{\alg})$. Each special quasi-endomorphism $x \in L(\calA_{\overline y})\otimes \bbQ$ gives rise to an element $ s(x )\in \omega_y$. Roughly speaking, $x$ induces an quasi-endomorphism $x[p^{\infty}]$ of the $p$-divisible group attached to $\calA_{\overline y}$, which, via Grothendieck-Messing theory, induces in turn a quasi-endomorphism $\xi = \xi(x)$ of the algebraic de Rham homology group $H_1^{dR}(\calA_{y})$ of the lift to characteristic zero. By virtue of being special,  $\xi$ lands in a distinguished quadratic subspace of quasi-endomorphisms that is identified with $V_{dR,y}$; the section $s(x)$ is obtained by pairing elements of $ \Fil^1  V_{dR, y} $ against $\xi$. All this is explained in more detail in \cite[\S 4.3]{AGHMP}, see also \cite{MadapusiPera}.

				Finally, by \cite[Proposition 4.3.2]{AGHMP}, if $x \in L(A_{\overline y})$ then the locus on which $x$ deforms to a special endomorphism of $\calA_{\tilde y}$ is precisely the vanishing locus $Z(s(x)) $ of $s(x)$; in particular,  applying \Cref{prop:Green Arakelov} and unwinding definitions, we have
				\begin{equation}
				\widehat \deg Z(\alpha,  \tau) =  \deghat  \mathscr Z(\alpha) \, q^{\alpha}
				\end{equation}
				for any totally positive $\alpha \in F^{\times}$ such that $\mathrm{Diff}(\alpha) = \{ v \}$ for a good place $v$.

				On the automorphic side, we need to modify the construction of our Eisenstein series slightly to account for the bad places. Define a section $\widetilde \Phi(g,s) = \otimes_{v \leq \infty} \widetilde \Phi_v(g,s) \in I(s, \chi)$ as follows:
				\begin{itemize}
					\item If $v$ is archimedean, we set $\widetilde \Phi_v(g,s) = \Phi^*_v(g,s)$ as in \Cref{sec:Phi*}.
					\item If $v$ is a non-archimedean place lying over the rational prime $p$, let $\varphi_{L^{\max},v}$ denote the characteristic function of $L^{\max}_p \cap \calV_v$, and let $\widetilde \Phi_v(g, s) = \Phi(g, s, \lambda_{\calV} (\varphi_{L^{\max},v}))$ denote the corresponding section. 
				\end{itemize}
				Note that for good places $v$, we have that $\widetilde \Phi_v(g,s) = \Phi^*_v(g,s)$, where $\Phi^*_v(g,s)$ is defined as in \Cref{sec:Phi*}.
				
				Then, after a mild modification of the proof of  \Cref{thm:main} to account for the bad places, we have the following result, which appears as Theorem 7.8.1 in \cite{AGHMP}, and is the main ingredient in the proof of the Bruinier-Kudla-Yang formula; to reiterate the essential point, our proof of this result does not require the explicit computation of the left hand side, and relies instead on the Siegel-Weil formula to establish the link to the right hand side.
				
				\begin{theorem}
					Let $\alpha \in F^{\times}$, and suppose that $\mathrm{Diff}(\alpha) = \{v\}$ for a good non-archimedean place $v$. Then
					\begin{equation}
					\deghat \mathscr Z(\alpha) q^{\alpha} = -  \deg_{\bbC}(Y_K) \, E_{\alpha}'(\tau, 0, \widetilde \Phi ).
					\end{equation}
					\qed
				\end{theorem}

\bibliographystyle{alpha}
\bibliography{refs}

\author{\noindent \small \textsc{Department of Mathematics, University of Manitoba} \\
	{\it E-mail address:} \texttt{siddarth.sankaran@umanitoba.ca}}

\end{document}